\documentclass[reqno]{amsart}
\usepackage{amssymb,color,epsf}
\usepackage [cmtip,arrow]{xy}
\xyoption{all}
\usepackage {pb-diagram,pb-xy}
\usepackage[mathscr]{eucal}

\ifx\pdfpageheight\undefined
   \usepackage[dvips,colorlinks=true,linkcolor=blue,citecolor=red,%
      urlcolor=green]{hyperref}
   \usepackage[dvips]{graphicx}
   \makeatletter
   \edef\Gin@extensions{\Gin@extensions,.mps}
   \makeatother
\else
 \usepackage[pdftex]{graphicx}
   \usepackage[bookmarksopen=false,pdftex=true,breaklinks=true,%
      backref=page,pagebackref=true,plainpages=false,%
      hyperindex=true,pdfstartview=FitH,colorlinks=true,%
      pdfpagelabels=true,colorlinks=true,linkcolor=blue,%
      citecolor=red,urlcolor=green,hypertexnames=false%
      ]%
   {hyperref}
\fi
\usepackage{bm}

\definecolor{DarkBlue}{rgb}{0,0.1,0.55}
\newcommand{\bfdef}[1]{{\bf \color{DarkBlue}  {#1}}}

\newtheorem{theorem}{Theorem}[section]

\newtheorem{corollary}[theorem]{Corollary}
\newtheorem{proposition}[theorem]{Proposition}
\theoremstyle{definition}
\newtheorem{definition}[theorem]{Definition}
\newtheorem{example}[theorem]{Example}

\newtheorem{question}[theorem]{Question}
\theoremstyle{remark}
\newtheorem{remark}[theorem]{Remark}
\numberwithin{equation}{section}

\newcommand {\sign} {\mbox{\bf sign}}

\newtheorem{notation}{\sc Notation}

\newcommand {\R} {\mbox{\rm R}}

\newcommand {\junk}[1]{}
\newcommand {\hide}[1]{}

\newcommand {\s}        {\mbox{\rm sign}}
\newcommand {\D}     {\mbox{\rm D}}

\newcommand {\C}     {\mathrm{C}}

\newcommand {\Real}[1]   {\mbox{${\mathbb R}^{#1}$}}

 \newcommand {\re}         {\Real{}}

\newcommand {\Z}  {{\mathbb Z}}
\newcommand {\Q}         {{\mathbb Q}}

\newcommand {\ZZ} {\mathrm{Zer}}
\newcommand {\RR} {\mathrm{Reali}}

\newcommand {\Der} {{\rm Der}}

\newcommand {\Ext} {{\rm Ext}}
\newcommand {\la}   {{\langle}}
\newcommand {\ra}   {{\rangle}}
\newcommand {\eps} {{\varepsilon}}
\newcommand {\E} {{\rm Ext}}

\def\sign{{\rm sign}}

\newcommand {\B} {\mbox {\rm B}}

\def\addots{\mathinner{\mkern1mu
\raise1pt\vbox{\kern7pt\hbox{.}}
\mkern2mu\raise4pt\hbox{.}\mkern2mu
\raise7pt\hbox{.}\mkern1mu}}

\newcommand {\Def} {\mbox {\rm Def}}

\newcommand{\RM}  {\mbox{\rm RM}}

\newcommand {\BElim} {\mbox {\rm BElim}}
\newcommand {\card} {\mbox {\rm card}}

\newcommand {\Pos} {\mbox {\rm Pos}}

\newcommand {\HH} {{\rm H}}
\newcommand {\LL} {{\rm L}}

\newcommand {\X} {\mathbf{X}}

\newcommand {\Sphere} {\mathbf{S}}

\newcommand{\coucou}[1]{\ifvmode\else\marginpar[\hfill$\rhd$]{$\lhd$}\fi
                       $\langle$\textsc{#1}$\rangle$}

\newcommand{\basu}[1]{}
\newcommand {\rp}[1]{}

\begin{document}
\title[]{
Algorithms in Real Algebraic Geometry: A Survey
\footnote{2000 Mathematics Subject Classification Primary 14P10, 14P25;
Secondary 68W30}
}
\author{Saugata Basu}
%%    Address of record for the research reported here
\address{Department of Mathematics,
Purdue University, West Lafayette, IN 47907, U.S.A.}
%%   Current address
\email{sbasu@math.purdue.edu}
%%\thanks{The author was supported in part by NSF grant CCF-0634907. 
%%\subjclass{Primary 14P10, 14P25; Secondary 68W30}
%%\date{January 1, 1994 and, in revised form, June 22, 1994.}

%%\dedicatory{This paper is dedicated to our advisors.}
\thanks{The author was  partially supported by an NSF grants CCF-0915954, CCF-1319080 and DMS-1161629.}
% ----------------------------------------------------------------

\begin{abstract}
We survey both old and new developments in the theory of algorithms
in real algebraic geometry -- starting from effective quantifier
elimination in the first order theory of reals due to Tarski and Seidenberg,
to more recent algorithms for computing topological invariants of
semi-algebraic sets. We emphasize throughout the complexity aspects
of these algorithms and also discuss the computational hardness of the
underlying problems. We also describe some recent results linking
the computational hardness of decision problems in the first order theory
of the reals, with that of computing certain topological invariants of 
semi-algebraic sets. Even though we mostly concentrate on 
exact algorithms, we also discuss some numerical 
approaches involving semi-definite programming that have gained popularity
in recent times.
\end{abstract}
\maketitle

% ----------------------------------------------------------------
\tableofcontents
\section{Introduction}
\label{sec:intro}
We survey developments in the theory of algorithms in real
algebraic geometry -- starting from the first effective quantifier elimination
procedure due to Tarski and Seidenberg, to more recent work
on efficient algorithms for quantifier elimination, as well as
algorithms for computing topological invariants of semi-algebraic sets -- such
as the number semi-algebraically connected components, Euler-Poincar\'e
characteristic, Betti numbers etc. Throughout the survey, the emphasis is on
the worst-case complexity bounds of these algorithms, and the continuing
effort to design algorithms with better complexity. Our goal in this survey 
is to describe these algorithmic results 
(including stating precise complexity bounds in  most cases), 
and also give some indications of the techniques involved in designing these
algorithms. 
We also 
describe some hardness results which show the intrinsic difficulty 
of some of these problems.

\subsection{Notation}
\label{subsec:notation}
We first fix some notation.
Throughout,  $\R$ will denote 
a \bfdef{real closed field} (for example, the field $\re$ of real numbers
or $\re_{{\rm alg}}$ of real algebraic numbers), and we will denote
by $\C$ the algebraic closure of $\R$.

A \bfdef{semi-algebraic subset} of $\R^k$ is a set defined by a
finite system of polynomial equalities and inequalities,
or more generally by a Boolean
formula whose atoms are polynomial equalities and inequalities.
Given a finite set ${\mathcal P}$ of polynomials in ${\R}[X_1,\ldots,X_k]$, a
subset $S$ of ${\R}^k$ is 
\bfdef{${\mathcal P}$-semi-algebraic}
if $S$ is the realization of a Boolean formula with
atoms $P=0$, $P > 0$ or $P<0$ with $P \in {\mathcal P}$
(we will call such a formula a quantifier-free $\mathcal{P}$-formula).

It is clear that for every semi-algebraic subset $S$ of $\R^k$
there exists
a finite set ${\mathcal P}$ of polynomials in ${\R}[X_1,\ldots,X_k]$
such that $S$ is ${\mathcal P}$-semi-algebraic.
We call a semi-algebraic set  a 
\bfdef{${\mathcal P}$-closed} semi-algebraic set
if it is defined by a Boolean formula with no negations with atoms
$P=0$, $P \geq 0$, or $P \leq 0$ with $P \in {\mathcal P}$.

For an element $a \in\R$ we let
$$
\s(a) =
\begin{cases}
0 & \mbox{ if }  a=0,\cr
1 & \mbox{ if } a> 0,\cr
-1& \mbox{ if } a< 0.
\end{cases}
$$

A  \bfdef{sign condition}  on
${\mathcal P}$ is an element of $\{0,1,- 1\}^{\mathcal P}$.
For any semi-algebraic set $Z \subset \R^k$ 
the \bfdef{realization of the sign condition $\sigma$ over $Z$}, 
$\RR(\sigma,Z)$, is the semi-algebraic set
\label{def:R(Z)}
$$
        \{x\in Z\;\mid\; \bigwedge_{P\in{\mathcal P}} \s({P}(x))=\sigma(P) \},
$$
and in case $Z= \R^k$ we will denote $\RR(\sigma,Z)$ by just $\RR(\sigma)$.
   
If ${\mathcal P}$ is a finite subset of
$\R [X_1, \ldots , X_k]$, we write the set of zeros
of ${\mathcal P}$ in $\R^k$ as
$$
\ZZ({\mathcal P},\R^k)=\{x\in \R^k\mid\bigwedge_{P\in{\mathcal P}}P(x)= 0\}.
$$

Given a semi-algebraic set $S \subset \R^k$, we will denote by
\bfdef{$b_i(S)$} the $i$-th \bfdef{Betti number} of $S$, 
that is the rank of the $i$-th homology group 
of $S$ (see \cite{BPRbook2}
for precise definitions of homology groups for semi-algebraic sets
defined over arbitrary real closed fields). 
Note that $b_0(S)$ is the number
of semi-algebraically connected components of $S$.
We will denote by 
$b(S)$ the sum $\sum_{i \geq 0} b_i(S)$.

For $x \in \R^k$ and $r > 0$, we will denote by $\B_k(x,r)$ (resp. 
$\Sphere^{k-1}(x,r)$) the open ball (resp. the sphere) with center $x$ and
radius $r$ in $\R^k$. 
%%When $x = 0$, we will write $\B_k(r)$ (resp. 
%%$\Sphere^{k-1}(r)$) instead of $\B_k(0,r)$ (resp. $\Sphere^{k-1}(0,r)$).
We will also denote the unit ball (resp. sphere) in $\R^k$ 
centered at $0$ by $\B_k$ (resp. $\Sphere^{k-1}$).

\subsection{Main algorithmic problems}
\label{subsec:problems}
Algorithmic problems in semi-algebraic geometry typically consist of 
the following. We are given as input a finite family, ${\mathcal P} \subset
\D[X_1,\ldots,X_k]$, where $\D$ is an ordered domain contained in the
real closed field $\R$.
The main algorithmic problems can be roughly divided into two classes
(though we will see later in Section \ref{sec:Toda} how they are 
related from the point of computational complexity).

The first class of problems has a logical flavor. It includes the 
following.

Given a quantified $\mathcal{P}$-formula
$\Phi$ (with or without free variables), the task is to:

\begin{enumerate}
\item (\bfdef{The Quantifier Elimination Problem})  Compute a quantifier-free formula
equivalent to $\Phi$.
\item (\bfdef{The General Decision Problem}) This is a special case of the previous
problem when $\Phi$ has no free variables, and the problem is to decide
the truth or falsity of $\Phi$.
\item
(\bfdef{The Existential Problem}) 
This is a special case of the last problem when there
is exactly one block of existential quantifiers; equivalently, the
problem can be stated as deciding whether a given $\mathcal{P}$-semi-algebraic
set is empty or not.
\end{enumerate}

The second class of problems has 
a more  geometric or topological flavor.
Given a description of a $\mathcal{P}$-semi-algebraic set $S \subset \R^k$
the task is to
decide whether certain geometric or  topological properties hold for $S$,
and in some cases also computing certain topological invariants of $S$.
Some of the most basic problems include the following. 

\begin{enumerate}
\item 
(\bfdef{Deciding Emptiness})
Decide whether $S$ is empty or not (this is the same as the
Existential Problem described above).
\item
(\bfdef{Deciding the existence of semi-algebraic connecting paths})
Given two points $x,y \in S$,  decide if they  belong to the
same semi-algebraically connected component of $S$.
\item
(\bfdef{Computing descriptions of semi-algebraic paths})
Given two points $x,y \in S$,  which  belong to the
same semi-algebraically connected component of $S$, output a description of a
semi-algebraic path in $S$ connecting $x,y$.
\item
\bfdef{(Describing Connected Components)}
Compute semi-algebraic descriptions of the semi-algebraically 
connected components of $S$.
\end{enumerate}

At a slightly deeper level we have  problems of a more topological flavor, such as:
\begin{itemize}
\item[(4)]
\bfdef{(Computing Betti Numbers)}
Compute  the cohomology groups of $S$,
its Betti numbers, its Euler-Poincar\'e characteristic  etc..
\item[(5)] 
\bfdef{(Computing Triangulations)}
Compute a semi-algebraic triangulation of $S$ 
as well as, 
\item[(6)]
\bfdef{(Computing Regular Stratifications)}
compute   a  partition of $S$ into smooth semi-algebraic subsets of various
dimensions satisfying certain extra regularity conditions
(for example, Whitney conditions (a) and(b)).
\end{itemize}

\begin{definition}[Complexity]
\label{def:complexity}
A typical  input to the algorithms considered in this survey 
will be a set  of polynomials with coefficients in an ordered 
domain $\D$ (which can be taken to be the ring generated by the
coefficients of the input polynomials).
By \bfdef{complexity of an algorithm} we will mean the function 
whose argument is the size of the input to the algorithm, 
measured by the number of number of polynomials, their degrees, and
the number of variables, and whose value is the supremum 
over all inputs of size equal to the argument, of the number of arithmetic
operations (including comparisons) performed by the algorithm
in the domain $\D$. 
In case the input polynomials have integer coefficients with bounded bit-size,
then we will often give the bit-complexity, which gives an upper bound on the number of bit
operations performed by the algorithm.
We refer the reader to \cite[Chapter 8]{BPRbook2} for a full discussion
about the various measures of complexity.
\end{definition}

The complexity of an algorithm  (see Definition \ref{def:complexity} 
above) for solving any of the
above problems is measured in terms of the following
three parameters:
\begin{itemize}
\item
the number of polynomials, $s = \card\;{\mathcal P}$,
\item
the maximum degree, $d = \max_{P \in {\mathcal P}} \deg(P)$, and 
\item
the number of variables, $k$ (and in case of quantifier elimination
problems, the block decomposition of the $k$ variables).
\end{itemize}

The rest of the paper is organized as follows. In Section \ref{sec:qe}, we
describe known algorithms for quantifier elimination in the first order theory of the
reals, starting from Tarski's algorithm, algorithms via cylindrical
algebraic decomposition, and finally more modern algorithms using the
critical points method. 
We discuss some variants of the quantifier elimination problem that arise
in applications, as well as
certain approaches using complex geometry
of polar varieties that give efficient probabilistic algorithms. 
We also discuss the known lower bounds for real  quantifier elimination.

In Section \ref{sec:top}, we concentrate on algorithms for computing
topological properties of semi-algebraic sets -- including connectivity
property via construction of roadmaps, computing the generalized 
Euler-Poincar\'e characteristic of semi-algebraic sets, as well as 
computing the Betti numbers of semi-algebraic sets. Throughout this section
the emphasis is on algorithms with singly exponential complexity bounds.
We also discuss certain results that are special to semi-algebraic sets
defined by quadratic inequalities, or more generally where the defining
polynomials have at most quadratic dependence on most of the variables.
We also point out the significance of some of the results from the point
of view of computational complexity theory. Finally, we discuss a 
recent reduction
result linking the complexity of the problem of computing the Betti numbers
of semi-algebraic sets, with that of the decision problem in the first order
theory of the reals with a fixed number of quantifier blocks.

In Section \ref{sec:sos}, we discuss numerical algorithms for polynomial
optimization using the ``sums-of-square'' approach. The main algorithmic
tool here is ``interior-point algorithms for semi-definite programming'' 
and we discuss the known results
on the computational complexity of the semi-definite programming problem.

We end with a list of open problems (Section \ref{sec:open}). \\

\noindent{\em Warning.}
There are several interesting topics which come under the purview
of algorithms in real algebraic geometry that have been left out of this
survey (because of lack of space as well as the author's lack of
expertise in some of these topics). 
For example, we do not make any attempt to survey the extremely
broad area of research concerning efficient implementation of theoretically
efficient algorithms, specific low dimensional applications such as computing 
the topology of curves and surfaces, computing certificates of positivity
of polynomials (for archimedean as well as non-archimedean real closed fields),
homotopy continuation algorithms for solving real systems etc. There are
multiple excellent sources available for most of these topics. 
Finally, algorithmic real algebraic geometry has a great variety of  
applications,
due to the ubiquity of semi-algebraic sets arising in different areas of
science and engineering -- including robotics, molecular chemistry, theoretical
computer science, database theory etc. We do not make any attempt to survey 
these applications.

\section{Quantifier elimination and related problems}
\label{sec:qe}
We begin appropriately with the first algorithm (in the modern sense)
in real algebraic geometry which is a starting point of the subject.
 
\subsection{The Tarski-Seidenberg Theorem and effective quantifier elimination}
Let 
$
\mathcal{P} = \{P_1, \ldots,P_s \} \subset
\R[X_1,\ldots,X_k,Y_1,\ldots,Y_{\ell}],
$
and $\Phi(Y)$ a first-order formula given by
\[
(Q_\omega X^{[\omega]}) \ldots (Q_{1}X^{[1]})F(P_1,\ldots,P_s),
\]
where
$Q_i \in \{ \forall,\exists \}, \; Q_i \neq Q_{i+1}$, 
$ Y = (Y_1,\ldots, Y_{\ell})$ is a block of $\ell$ free variables, 
$X^{[i]}$ is a block of ${k_i}$ variables with 
$\sum_{1\leq i\leq \omega}k_i =k,$ and
$F(P_1,\ldots,P_s)$ is a quantifier-free Boolean formula with atomic
predicates of the form 
$\sign(P_i (Y,X^{[\omega]}, \ldots , X^{[1]})) = \sigma $ 
where $\sigma \in \{ 0, 1 , -1 \}$. 
If $\Pi$ denotes the partition
of the blocks of variables $X_1,\ldots,X_k$ into 
the $\omega$ blocks of sizes $k_1,\ldots,k_\omega$,
in the formula $\Phi$, we call $\Phi$ 
a \emph{$(\mathcal{P},\Pi)$-formula}.

The Tarski-Seidenberg theorem states that

\begin{theorem}
\cite{Tarski51}
There exists 
a quantifier-free
formula,
$\Psi(Y)$, such that for any ${y} \in R^\ell,$ $\Phi(y)$ is true if and
only if $\Psi({y})$ is
true. 
\end{theorem}

The quantifier elimination problem is to algorithmically construct
such a formula. 

\subsubsection{Algorithm arising from Tarski's proof}
Tarski's proof \cite{Tarski51}  
of the existence of quantifier elimination in the first order
theory of the reals is effective and is based on Sturm's theorem
for counting real roots of polynomials in one variable with real
coefficients used in a parametric way. A modern treatment of this
proof can be found in \cite[Chapter 2]{BPRbook2}. The complexity
of this procedure was not formally analysed in Tarski's paper. However, 
the algorithm eliminates one variable at a time using a parametrized 
version of Euclidean remainder sequence, and as a result the number and
degrees of the polynomials in the remaining variables grow rather fast,
and it is not possible to bound the complexity of the algorithm by
any function which is a tower of exponents (in the input parameters) 
of a fixed height, which implies that the complexity of Tarski's
algorithm is not {\em elementary recursive}. An elementary recursive
algorithm for the General Decision Problem was found later 
by Monk \cite{monk}.

\subsubsection{Cylindrical Algebraic Decomposition}
\label{sec:cad}
One fundamental technique for computing topological
invariants of semi-algebraic sets is  \emph{Cylindrical Algebraic
Decomposition}. 
Even though the mathematical ideas behind cylindrical algebraic decomposition
were  known before (see for example \cite{Loj2}),
Collins and W\"uthrich  \cite{Col, Collins75,Wuthrich76} were the first to
apply cylindrical algebraic decomposition in the setting of 
algorithmic semi-algebraic geometry.
Schwartz and Sharir \cite{SS}
realized its importance in trying to solve the motion planning problem in 
robotics, as well as computing topological properties of semi-algebraic
sets. 

\begin{definition}[Cylindrical Algebraic Decomposition]
\label{5:def:cad}
A \bfdef{cylindrical algebraic decomposition} of ${\R}^k$
is a sequence ${\mathcal S}_1,\ldots,{\mathcal S}_k$
where, for each $1\leq i\leq k$, ${\mathcal S}_i$ is a finite
partition of ${\R}^i$ into  semi-algebraic
subsets, called the  cells of level $i$,
which satisfy the following properties:
\begin{itemize}
\item Each cell $S\in {\mathcal S}_1$ is either a point or
 an open interval.
\item For every $1\leq i<k$ and every $S\in {\mathcal S}_i$,
there are finitely
 many continuous
semi-algebraic functions
 $$ \xi_{S,1}<\ldots<\xi_{S,\ell_S}: S\longrightarrow {\R}$$ such that
the cylinder $S\times {\R} \subset {\R}^{i+1}$ is the disjoint union
 of cells of ${\mathcal S}_{i+1}$ each of 
which is:
\begin{itemize}
\item either the graph
of one of the functions $\xi_{S,j}$, for $j =
1,\ldots,\ell_S$: $$
\{ (x', x_{j+1}) \in S \times {\R}\mid x_{j+1} = \xi_{S, j}(x') \}\;,
$$
\item or a  band
 of the cylinder bounded from below and from above by the graphs
 of the functions
$\xi_{S,j}$ and $\xi_{S,j+1}$, for $j = 0,\ldots, \ell_S$,
 where we take $\xi_{S,0} =-\infty$
and $\xi_{i,\ell_S+1} = +\infty$:
$$
\{ (x', x_{j+1}) \in S \times {\R}\mid \xi_{S,j}(x')
< x_{j+1} < \xi_{S, j+1}(x')\}\;.
$$
\end{itemize}
\end{itemize}
Given a semi-algebraic subset $S \subset \R^k$, we say that a cylindrical algebraic decomposition, 
$\mathcal{S}_1,\ldots,\mathcal{S}_k$, of $\R^k$, is \emph{adapted to $S$}, if for every cell $C \in \mathcal{S}_k$, either $C \subset S$ or $C \cap S = \emptyset$.
\end{definition}

\begin{definition}
\label{def:invariantcad}
Given a finite set ${\mathcal P} \subset {\R}[X_1,\ldots,X_k]$, a
subset $S$ of ${\R}^k$ is 
is \bfdef{${\mathcal P}$-invariant} 
if every polynomial $P \in {\mathcal P}$
has a constant sign ($>0$, $<0$, or $=0$) on $S$.
A \bfdef{cylindrical algebraic decomposition of ${\R}^k$ adapted to
\index{Cylindrical decomposition!adapted to ${\mathcal P}$}
\label{5:de:cylindricaladap}
${\mathcal P}$} is a cylindrical algebraic decomposition for which each cell
 $C\in {\mathcal S}_k$ is ${\mathcal P}$-invariant.  It is clear that if
$S$ is
${\mathcal P}$-semi-algebraic, a cylindrical algebraic 
 decomposition adapted to
${\mathcal P}$
is a cylindrical algebraic 
 decomposition adapted to
$S$.
\end{definition}

One important result which underlies most algorithmic applications of 
cylindrical algebraic decomposition is the following
(see \cite[Chapter 11]{BPRbook2} for an easily accessible exposition).

\begin{theorem} \cite{Collins75, Wuthrich76}
\label{5:the:cad} 
For every finite set ${\mathcal P}$
of polynomials in ${\R}[X_1,\ldots,X_k]$, there is a cylindrical
decomposition of ${\R}^k$ adapted to ${\mathcal P}$.
Moreover, such a decomposition can be computed in time
$(sd)^{2^{O(k)}}$, where $s = \card \;{\mathcal P}$ and 
$d = \max_{P \in {\mathcal P}} \deg(P)$.
\end{theorem}

Cylindrical algebraic decomposition provides an alternative (and more
efficient compared to Tarski's) 
algorithm for quantifier elimination, since (using the same
notation as in the previous section) the semi-algebraic subset 
of $\R^\ell$ defined by $\Phi(Y)$, is a union of cells (of various
dimensions) in a
cylindrical algebraic decomposition of $\R^{k+\ell}$ 
adapted to $\mathcal{P}$ (cf. Definition \ref{def:invariantcad}),
where $Y_1,\ldots,Y_\ell$ are the last $\ell$ variables. 
This last fact is a consequence of the 
``cylindrical'' structure of the decomposition. The complexity of such
an algorithm is bounded by the complexity of
computing the cylindrical decomposition and is doubly exponential.
More precisely, the complexity is bounded by $(s d)^{2^{O(k+\ell)}}$.

\begin{remark}
The technique of cylindrical algebraic decomposition is also used
in algorithms for computing topological properties of semi-algebraic sets.
After making a generic linear change of co-ordinates, the cylindrical
algebraic decomposition algorithm yields a finite cell complex from which
topological invariants of the underlying semi-algebraic sets can be extracted.
It should be noted that a change of co-ordinates is needed to 
obtain a cell complex. However, in certain applications 
a change of co-ordinates is not allowed  (see \cite{BGV2013} for 
one such application).
\end{remark}

\begin{question}
\label{ques:regular-cad}
Does there always exists a semi-algebraic 
\emph{cell decomposition}  adapted to a given finite family of 
polynomials, having a cylindrical structure with respect to the 
any given set of coordinates?
\end{question}

\subsubsection{Lower bound}
Given the doubly exponential upper bound on the complexity of
quantifier elimination algorithm that follows from cylindrical
algebraic decomposition, it is interesting to ask 
whether it is at all possible to do better.
This question was 
investigated by Davenport and Heintz \cite{Davenport-Heintz88} 
who proved a doubly exponential {\em lower bound}  on the complexity of real
quantifier elimination, by constructing a sequence of quantified
formula having the property that any equivalent sequence of 
quantifier-free formulas would necessarily have 
doubly exponential growth in size. However, the quantified formulas
in the sequence they constructed had a large number of quantifier
alternations (linear in the number of variables). Thus, 
while it is impossible to hope for better than doubly exponential
dependence in the number, $\omega$,  of quantifier alternations, 
it might still be possible to obtain algorithms with much better
complexity (i.e. singly exponential in the number of variables)  
if we fix the number of quantifier alternations. This is what
we describe next.

\subsection{The critical points method and singly exponential algorithms}
\label{sec:critical}
As mentioned earlier, all algorithms using cylindrical algebraic decomposition
have doubly exponential complexity. Algorithms with singly exponential
complexity for solving problems in semi-algebraic geometry are 
mostly based on the {\em critical points method}. This method was
pioneered by Grigoriev and Vorobjov and 
\cite{GV,GV92},
Renegar \cite{R92}, and later improved in various ways by several researchers including Canny \cite{Canny93a}, 
Heintz, Roy and Solern\`o \cite{HRS94},
Basu, Pollack and Roy 
\cite{BPR95} amongst others.
In simple terms, the \bfdef{critical points method} is nothing but
a method for finding at least one point
in every semi-algebraically  connected
component of an algebraic set.
It can be shown that for a bounded nonsingular algebraic hyper-surface, 
it is possible to change coordinates
so that its projection to the $X_1$-axis has
a finite number of non-degenerate critical points.
These points
provide at least one point in every semi-algebraically
 connected component of the
bounded nonsingular algebraic
hyper-surface.
Unfortunately this is not very useful
in algorithms since it provides no method
for performing this linear change of variables. Moreover
when we deal with the case of a
general algebraic set, which may be unbounded or singular,
this method no longer works.

In order to reduce the general case to the case of bounded nonsingular
algebraic sets, we use  
an important technique in algorithmic semi-algebraic geometry --
namely, perturbation of a given real algebraic set in $\R^k$ using 
one or more infinitesimals. 
The perturbed variety is then defined over a non-archimedean real
closed extension of the ground field -- namely the field of algebraic
Puiseux series in the infinitesimal elements with coefficients in $\R$.

Since the theory behind such extensions might be unfamiliar to some
readers, we introduce here the necessary algebraic background 
referring the reader to \cite[Section 2.6]{BPRbook2} 
for full detail and proofs.

\subsubsection{Infinitesimals and the field of algebraic Puiseux series}
\label{subsec:Puiseux}

\begin{definition}[Puiseux series]
A \bfdef{Puiseux series} in $\eps$ with coefficients in $\R$ 
is a  series of the form
\begin{equation}
\label{2:not:puiseux1}
		{\overline a}=\sum_{i \ge k} a_i \eps^{i/q},
\end{equation}
with $k \in\Z$, $i\in\Z$, $a_i\in \R$, $q$ a positive integer.
\end{definition}

It is a straightforward exercise to verify that
the field of all Puiseux series in $\eps$ with coefficients in $\R$ 
is an ordered field, in which the order extends the
order of $\R$, and $\eps$ is an infinitesimally small and positive element
i.e.,  is positive and smaller than any positive $r\in \R$.

\begin{notation}
\label{not:puiseux}
The field of Puiseux series in $\eps$ with coefficients in $\R$ contains
as a subfield, the field of Puiseux series which are algebraic over
$\R[\eps]$.
We denote by $\R\langle \eps\rangle$  the \bfdef{field of algebraic
Puiseux series} in $\eps$ with coefficients in $\R$.
We will also use the notation $\R\la\eps_1,\ldots,\eps_m\ra$ to denote
the field $\R\la\eps_1\ra\cdots\la\eps_m\ra$. Notice that in the
field $\R\la\eps_1,\ldots,\eps_m\ra$
for $1\leq i \leq m$, 
the element $\eps_i \in \R\la\eps_1,\ldots,\eps_i\ra \subset \R\la\eps_1,\ldots,\eps_m\ra$ 
is positive and infinitesimal over the subfield $\R\la\eps_1,\ldots,\eps_{i-1}\ra$ of 
$\R\la\eps_1,\ldots,\eps_{i}\ra$.
 \end{notation}

The following theorem is classical (see for example
\cite[Section 2.6]{BPRbook2} for a proof).

\begin{theorem}
\label{the:Puiseux}
The field $\R\langle \eps\rangle$ is real closed.
\end{theorem}

\begin{definition}[The $\lim_\eps$ map]
\label{def:lim}
When $a \in \R\la \eps \ra$ is bounded from above and below by an element of $\R$,
$\lim_\eps(a)$ is the constant term of $a$, obtained by
substituting 0 for $\eps$ in $a$.
\end{definition}

\begin{example}
\label{ex:lim}
A typical example of the application of the $\lim$ map can be
seen in Figures \ref{fig:homotopy1} and  \ref{fig:homotopy2} below. 
The first picture depicts the algebraic set $\ZZ(Q,\R^3)$.
The second picture depicts 
the algebraic set $\ZZ(\Def(Q,\zeta,4),\R\la\zeta\ra^3)$ (where we substituted
a very small positive number for $\zeta$ in order to be able to display
this set).
The polynomials $Q$ and $\Def(Q,\zeta,4)$ 
are defined by \eqref{eqn:examplehomotopy} and
\eqref{11:equation:deform1} respectively.
The algebraic sets $\ZZ(Q,\R^3)$ and $\ZZ(\Def(Q,\zeta,4),\R\la\zeta\ra^3)$
are related by
\[
\ZZ(Q,\R^3) = \lim_\zeta~\ZZ(\Def(Q,\zeta,4),\R\la\zeta\ra^3).
\]
\end{example}

Since we will often consider the semi-algebraic sets
defined by the same formula, but over different real closed
extensions of the ground field, the following notation is 
useful.

\begin{notation}
\label{not:extension}
Let $\R'$ be a real closed field containing $\R$.
Given a semi-algebraic subset
$S \subset \R^k$, defined by a first-order formula $\Phi$ (with coefficients in $\R$),  
the \bfdef{extension}
of $S$ to $\R'$, denoted $\E(S,\R')$, is
the semi-algebraic subset of $\R'^k$ defined by the same
formula $\Phi$.
\end{notation}
The set $\E(S,\R')$ is well defined (i.e. it only depends on the set
$S$ and not on the quantifier free formula chosen to describe it).
This is an easy consequence of the Tarski-Seidenberg transfer principle 
(see for example \cite[Theorem 2.80]{BPRbook2}).

We now return to the discussion of the critical points method.
In order for the critical points method to work for all algebraic sets,
we associate to a possibly unbounded
algebraic set $Z \subset \R^k$ a bounded algebraic set
$Z_{\mathrm{b}} \subset \R\la \eps \ra^{k+1},$
whose semi-algebraically connected components are closely related to those
of $Z$.

Let $Z=\ZZ(Q,\R^k)$ and consider
$$Z_{\mathrm{b}}= \ZZ(Q^2+(\eps^2(X_1^2+\ldots+X_{k+1}^2)-1)^2,\R\la \eps \ra^{k+1}).$$
The variety  
$Z_{\mathrm{b}}$ is the intersection of the  sphere $S^{k}(0,1/\eps)$ of
center $0$ and radius $\displaystyle{1 \over \eps}$
with a cylinder
based $\Ext(Z,\R \la \eps \ra)$ (and is hence
bounded over $\R \la \eps \ra$).
The intersection of $Z_{\mathrm{b}}$ with the hyper-plane $X_{k+1}=0$
is the intersection  of $Z$ with  the  sphere $S^{k-1}(0,1/\eps)$ of
center $0$ and radius $\displaystyle{1 \over \eps}$.
Denote by $\pi$ the projection from
$\R\la \eps \ra^{k+1}$ to $\R\la \eps \ra^{k}.$

The following proposition which appears in \cite{BPRbook2} 
relates the semi-algebraically connected components of $Z$ with
those of $Z_{\mathrm{b}}$. This allows us to reduce the problem of finding
points on every semi-algebraically connected component of a possibly unbounded algebraic set
to the same problem on bounded algebraic sets.

\begin{proposition}
\label{11:prop:boundedunbounded}
Let $N$ be a finite number of points meeting  every semi-algebraic\-ally
connected component of $Z_{\mathrm{b}}$.
Then $\pi(N)$ has a non-empty intersection with  every semi-algebraically
connected component of the extension $\E(Z,\R\la \eps \ra)$.
\end{proposition}

We obtain immediately
using Proposition \ref{11:prop:boundedunbounded}
a method for finding a point in every semi-algebraically connected component
of an algebraic set. Note  that these points
have coordinates in the extension $\R \la \eps\ra$
rather than in the real closed field $\R$ we started with.
However, the extension from $\R$ to $\R \la \eps\ra$
preserves  semi-algebraically connected components.

\subsubsection{Representation of points}
One important aspect in any algorithm in real algebraic geometry
is how to represent points whose co-ordinates belong to some
real algebraic extension of the ordered ring $\D$ generated by
the coefficients of the input polynomials. There are as usual
several options, such as representing an arbitrary real algebraic
number using isolating intervals, or by Thom encodings etc.
In the singly-exponential algorithms described in the book \cite{BPRbook2},
points in $\R^k$ are represented by {\em univariate representations}
and an associated {\em Thom encoding}. 
Even though we will not need any further detail about these representations
in this survey, given their importance in most of the
algorithms that we refer to, 
we include their precise definitions below.

\begin{definition}[Thom encoding]
Let $P \in \R[X]$ and $\sigma \in \{0,1,-1\}^{\Der(P)}$,
a sign condition on the
 set $\Der(P)$
of derivatives of $P$.
The sign condition $\sigma$ is
a \bfdef{Thom encoding}  of
${x \in \R}$ if $\sigma(P)=0$ and
 $\sigma$  is the sign condition taken by
the set $\Der(P)$ at $x$.
Given a Thom encoding $\sigma$, we  denote
by  $x(\sigma)$ the root of $P$ in $\R$ specified by $\sigma$.
\end{definition}

(Note that the use of Thom encoding to
represent algebraic numbers was introduced in
algorithmic real algebraic geometry by Coste and Roy in \cite{CR}.)

\begin{definition}
[Univariate representations and real univariate representations]
A \bfdef{$k$-univariate representation}
is a $(k+2)$-tuple of polynomials of
$\R[T]$,
$$
(f(T),g_0(T),g_1(T),\ldots ,g_k(T)),
$$
such that $f$ and $g_0$ are coprime.

The \bfdef{points associated}
to a univariate
representation are the points 
$$
\Bigg( {g_1(t) \over g_0(t)},\ldots, {g_k(t)
\over g_0(t)}\Bigg)\in \C^k 
$$ where
$t\in\C$ is a root of $f(T)$.

A \bfdef{real $k$-univariate representation}
is a pair $u,\sigma$ where
$u$ is a $k$-univariate representation and $\sigma$ is the Thom
encoding of a root of $f$,
$t_\sigma \in \R$.
The \bfdef{point associated} to the real
univariate representation is the point
$$
\Bigg( {g_1(t_\sigma) \over g_0(t_\sigma)},\ldots, {g_k(t_\sigma) \over
g_0(t_\sigma)}\Bigg)\in
\R^k.
$$
\end{definition}

\begin{remark}
By parametrizing the definition of a real $k$-univariate representation
(lets say by a co-ordinate function such as $X_1$) one obtains 
descriptions of semi-algebraic curves. These \bfdef{curve segment
representations} play an important role in algorithms for computing
roadmaps of semi-algebraic sets (see Section \ref{sec:roadmap} below).
\end{remark}

\subsubsection{Deformation techniques to deal with singular varieties}
For dealing with possibly singular algebraic sets
we define \bfdef{$X_1$-pseudo-critical points} of $\ZZ(Q,\R^k)$ when
$\ZZ(Q,\R^k)$ is a bounded algebraic set.
These pseudo-critical points are a finite set of
points meeting every semi-algebraically connected
 component of $\ZZ(Q,\R^k)$. They are the
limits of the  critical points of the projection to
the $X_1$ coordinate of  a bounded nonsingular
algebraic hyper-surface defined by a particular
infinitesimal perturbation, $\Def(Q,\zeta,d)$, 
of the polynomial $Q$ (where $d =\deg(Q)$). Moreover, the
equations defining the critical points of
the projection on the $X_1$ coordinate
on the perturbed algebraic set have a very  special
algebraic structure (they form a Gr\"{o}bner basis 
\cite[Section 12.1]{BPRbook2}), 
which makes possible efficient computation 
of these pseudo-critical values and points. We refer the reader to
\cite[Chapter 12]{BPRbook2} for a full exposition 
including the definition
and basic properties of Gr\"{o}bner basis.

The deformation $\Def(Q,\zeta,d)$ of $Q$ is defined as follows.
Suppose that $\ZZ(Q,\R^k)$ is contained in the ball of center
$0$ and radius $1/c$.
Let $\bar d$ be an even integer bigger than the degree $d$ of $Q$ and let
\begin{equation}
G_k(\bar d,c)
=c^{\bar d}(X_1^{\bar d}+\cdots+
X_{k}^{\bar d}+X_2^2+\cdots+X _k^2)-(2k-1),
\end{equation}
\begin{equation}
\label{11:equation:deform1}
\Def(Q,\zeta,d) =\zeta G_k(\bar d,c)+{(1-\zeta)
}Q.
\end{equation}

The algebraic set $\ZZ(\Def(Q,\zeta,d),\R \la \zeta \ra ^k)$
is a bounded and non-singular 
hyper-surface lying infinitesimally close to $\ZZ(Q,\R^k)$
and the critical points of the projection
map onto the $X_1$ co-ordinate restricted to
$\ZZ(\Def(Q,\zeta,d),\R \la \zeta \ra ^k)$
form a finite set of points.
We take the 
images of these points under $\lim_\zeta$ (cf. Definition
\ref{def:lim}) and
we call the points obtained in this manner  the $X_1$-pseudo-critical points
of $\ZZ(Q,\R^k)$.
Their projections on the $X_1$-axis  are called pseudo-critical values.

\begin{example}
\label{ex:homotopy}
We illustrate the perturbation mentioned above by a concrete example.
Let $k=3$ and $Q \in \R[X_1,X_2,X_3]$ be defined by
\begin{equation}
\label{eqn:examplehomotopy}
Q = X_2^2 - X_1^2 + X_1^4 + X_2^4 + X_3^4.
\end{equation}
Then, $\ZZ(Q,\R^3)$ is a bounded algebraic subset of $\R^3$ shown below
in Figure \ref{fig:homotopy1}.
Notice that $\ZZ(Q,\R^3)$ has a singularity at the origin.
The surface  $\ZZ(\Def(Q,\zeta,4),\R\la\zeta\ra^3)$ 
with a small positive real number substituted for $\zeta$ is shown in
Figure \ref{fig:homotopy2}. 
Notice that this surface is non-singular, but has a different 
semi-algebraic homotopy
type than $\ZZ(Q,\R\la\zeta\ra^3)$ 
(it has three semi-algebraically connected components compared to
only one of  $\ZZ(Q,\R\la\zeta\ra^3)$). However, the semi-algebraic set bounded
by $\ZZ(\Def(Q,\zeta,4),\R\la\zeta\ra^3)$  (i.e. the part inside the
larger component but outside the smaller ones)
is semi-algebraically homotopy equivalent to $\ZZ(Q,\R\la\zeta\ra^3)$.

\vspace*{2.75cm}
\begin{figure}[hbt]
\centerline{\scalebox{0.3}{
\begin{picture}(500,100)
\includegraphics[bb=0 0 64mm 30mm]{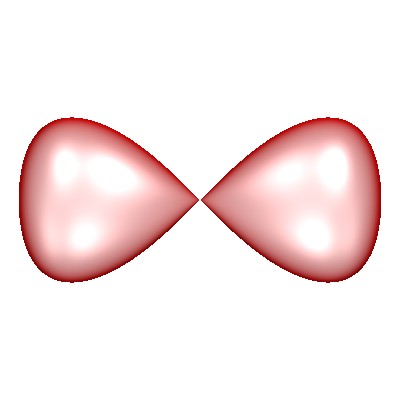}%
\end{picture}
}
}
\caption{The algebraic set $\ZZ(Q,\R^3)$.}
\label{fig:homotopy1}
\end{figure}

\vspace*{2.75cm}
\begin{figure}[hbt]
\centerline{\scalebox{0.3}{
\begin{picture}(500,100)
\includegraphics[bb=0 0 64mm 30mm]{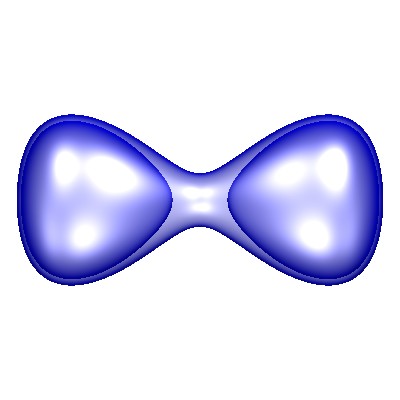}%
\end{picture}
}
}
\caption{The algebraic set $\ZZ(\Def(Q,\zeta,4),\R^3)$.}
\label{fig:homotopy2}
\end{figure}
\end{example}

By computing algebraic representations (see \cite[Section 12.4]{BPRbook2}
for the precise definition of such a representation)
of the pseudo-critical points
one obtains for any given algebraic set a finite set of points guaranteed to
meet every semi-algebraically connected component of this algebraic set. Using some 
more arguments from real algebraic geometry one can also reduce the problem
of computing a finite set of points guaranteed to meet every semi-algebraically connected
component of the realization of every realizable sign condition on a given
family of polynomials to finding points on certain algebraic sets defined by
the input  polynomials (or infinitesimal perturbations of these polynomials).
The details of this argument can be found in \cite[Proposition 13.2]
{BPRbook2}. 

The following theorem which is the best result of this kind appears in 
\cite{BPR95b}.

\begin{theorem}%%\cite{BPR95b}
\label{13:the:samplealg}
Let $\ZZ(Q,\R^k)$ be an algebraic set
of real dimension $k'$, where $Q \in \R[X_1, \ldots ,X_k]$, with $\deg(Q) \leq d$.  Let ${\mathcal P} \subset\R[X_1, \ldots ,X_k]$ be a set of $s$
polynomials with $\deg(P) \leq d$ for 
each $P\in {\mathcal P}$.
Let $\D$ be the
ring generated by the coefficients of $Q$ and the
 polynomials in ${\mathcal P}$.
There is an algorithm which
computes a set of points
meeting every
semi-algebraically connected
component of every  realizable sign condition on $\mathcal P$
over $\ZZ(Q,\R\la \eps,\delta \ra ^k)$. The algorithm has
complexity
$$(k'(k-k')+1)\displaystyle{\sum_{j \le k'}4^j {s  \choose j}
d^{O(k)}=s^{k'} d^{O(k)}}$$
 in $\D$.
There is also an
algorithm providing the
list of
signs of all the polynomials of
$\mathcal P$
at each of these points
 with
complexity
$$
(k'(k-k')+1)\displaystyle s{\sum_{j \le k'}4^j {s  \choose j}
d^{O(k)}=s^{k'+1} d^{O(k)}}
$$ in $\D$.
\end{theorem}

Notice that the combinatorial complexity 
(i.e. the part that depends on $s$) of the algorithm in Theorem
\ref{13:the:samplealg} depends on the dimension of the variety rather
than that of the ambient space. 

\subsection{Certain quantitative results on the topology of semi-algebraic sets}
\label{subsec:quantitative-topology}
Recall that for any semi-algebraic subset $S \subset \R^k$, we denote by $b_0(S)$ the number of semi-algebraically
connected components of $S$ (see Subsection \ref{subsec:notation}).

The problem of designing algorithms for computing sample points in every semi-algebraically connected
components of the realizations of sign conditions of a family of polynomials has a close connection
with the purely quantitative question of proving tight upper bounds on the number of such components
in terms of the number of polynomials in the family and their degrees. Such upper bounds serve as guiding
principles in design of algorithms, as the complexity of a sampling algorithm cannot be smaller than the number of such components in the worst case. Indeed the critical point method (along with certain tools
from algebraic topology such as the Mayer-Vietoris exact sequence) is a key element in the proofs of
such bounds. The following upper bound appears in \cite{BPR8}.

\begin{theorem}
\label{thm:bound-on-sign-conditions}
Let $\mathcal{P}, \mathcal{Q} \subset
  \R [X_1, \ldots, X_k]$ be finite families of
  polynomials such that the degrees of the polynomials in $\mathcal{P},
  \mathcal{Q}$ are bounded by $d$,
  $\ensuremath{\operatorname{card}}\mathcal{P}= s$, and
  $\dim_{\R} (V) = k'$, where $V
  =\ensuremath{\operatorname{Zer}} (\mathcal{Q},
  \R^k)$. Then,
\[
    \sum_{\sigma \in \{ 0, 1, - 1 \}^{\mathcal{P}}} b_0
    (\ensuremath{\operatorname{Reali}} (\sigma, V))  \leq  O (1)^k s^{k'} d^k .
\]
\end{theorem}

Notice that in Theorem \ref{thm:bound-on-sign-conditions} the bound depends on the maximum amongst the degrees of all the polynomials involved. It turns out that in recent applications, such as in 
discrete and computational geometry, a more refined dependence on the degrees is needed. In particular,
the situation where the degrees of the polynomials defining the variety $V$ is much smaller than the
degrees of the polynomials in $\mathcal{P}$ is of interest (see \cite{Solymosi-Tao,KMSS,Zahl}). The following theorem that achieves this is proved in \cite{Barone-Basu11a}.

\begin{theorem}
\label{thm:refined-bound-on-sign-conditions}
 Let $\mathcal{P}, \mathcal{Q} \subset
  \R [X_1, \ldots, X_k]$ be finite subsets of
  polynomials such that $\deg (Q) \leq d_1$ for all $Q \in \mathcal{Q}$, $\deg
  (P) \leq d_2$ for all $P \in \mathcal{P}$. Suppose also that $d_1 \leq
  d_2$, and the real dimension of $V =\ensuremath{\operatorname{Zer}}
  (\mathcal{Q}, \R^k)$ is $k_1 \leq k$, and that
  $\ensuremath{\operatorname{card}}\mathcal{P}= s$. Then,
  \begin{eqnarray*}
    \sum_{\sigma \in \{ 0, 1, - 1 \}^{\mathcal{P}}} b_0
    (\ensuremath{\operatorname{Reali}} (\sigma, V)) & \leq & O (1)^k (s
    d_2)^{k_1} d_1^{k - k_1} .
  \end{eqnarray*}
\end{theorem}

\begin{question}
\label{ques:refined_algo}
In view of Theorem \ref{thm:refined-bound-on-sign-conditions}, it is natural to ask whether there
is an algorithm, that given as input two finite sets $\mathcal{P}, \mathcal{Q} \subset
  \R [X_1, \ldots, X_k]$ 
 with $\deg (Q) \leq d_1$ for all $Q \in \mathcal{Q}$, $\deg
  (P) \leq d_2$ for all $P \in \mathcal{P}$, $d_1\leq d_2$,
computes a finite set points meeting every semi-algebraically connected component
of every  realizable sign condition on $\mathcal{P}$
over $\ZZ(\mathcal{Q},\R^k)$, and 
whose complexity has a more refined dependence on the degrees
$d_1$ and $d_2$ (reflecting the bound in Theorem \ref{thm:refined-bound-on-sign-conditions}) 
than what one obtains from Theorem \ref{13:the:samplealg}.
\end{question}

\subsection{Certain quantitative results in metric semi-algebraic geometry}
In the case $\D =\mathbb{Z}$, a careful analysis of the algorithm in
Theorem \ref{13:the:samplealg} produces an explicit upper bound on the
radius of a ball centered at the origin which is guaranteed to meet every
semi-algebraically connected component of any $\mathcal{P}$-semi-algebraic
set in terms of $s$, $d$, $k$ and a bound on the bit-size, $\tau$, of the
coefficients of $\mathcal{P}$. This and related bounds of this type 
are often needed
in designing other algorithms (for instance, in order to compute certificates
of positivity by sub-division method as done in \cite{BoCaR}). The following 
rather technical but completely explicit estimate appears
in \cite{BR10} (the same paper contains several other explicit estimates
of similar types).

\begin{notation}
Given an integer $n$, we denote by $\mathrm{bit} (n)$ the number of bits of its
absolute value in the binary representation. Note that
\begin{equation}
  \mathrm{bit} (nm) \leq \mathrm{bit} (n) + \mathrm{bit} (m), \label{bit1}
\end{equation}
\begin{equation}
  \mathrm{bit} \left( \sum_{i = 1}^n m_i \right) \leq  \mathrm{bit} (n) +
  \sup^n_{i = 1} (\mathrm{bit} (m_i)). 
\end{equation}
\end{notation}

\begin{theorem}\cite{BR10}
 \label{13:the:saball}
Let $\mathcal{P}=\{P_1, \ldots, P_s \} \subset
  \mathbb{Z}[X_1, \ldots, X_k]$ and suppose that $P \in \mathcal{P}$ have
  degrees at most $d$, and the coefficients of $P \in \mathcal{P}$ have
  bitsizes at most $\tau$. 
Then there exists a ball centered at the origin of
  radius
  \[ \left( (2 D N (2 N - 1) + 1) 2^{(2 N - 1) (\tau'' + \mathrm{bit} (2 N - 1)
     + \mathrm{bit} (2 DN + 1))} \right)^{1 / 2} \]
  where
  \begin{eqnarray*}
    d' & = & \sup (2 (d + 1), 6),\\
    D & = & k (d' - 2) + 2,\\
    N & = & d' (d' - 1)^{k - 1},\\
    \tau'' & = & N (\tau_2' + \mathrm{bit} (N) + 2 \mathrm{bit} (2 D + 1) + 1),\\
    \tau_2' & = & \tau_1' + 2 (k - 1) \mathrm{bit} (N) + (2 k - 1) \mathrm{bit}
    (k),\\
    \tau_1' & = & D (\tau_0' + 4 \mathrm{bit} (2 D + 1) + \mathrm{bit} (N)) - 2
    \mathrm{bit} (2 D + 1) - \mathrm{bit} (N),\\
    \tau_0' & = & 2 \tau + k \mathrm{bit} (d + 1) + \mathrm{bit} (2 d') +
    \mathrm{bit} (s)
  \end{eqnarray*}
  intersecting every semi-algebraically connected component of the realization
  of every realizable sign condition (resp. realizable weak sign condition)
  on  $\mathcal{P}$.
\end{theorem}

\begin{remark}
  Note that asymptotic bounds of the form $2^{\tau d^{O (k)}}$ 
  for the same problem were known before \cite{BPRbook2,GV,R92}. One 
  point which needs some explanation is the fact that $s$ plays a role in
  the estimate in Theorem \ref{13:the:saball}, while it does not appear in the
  formula \ $2^{\tau d^{O (k)}}$. This is because the total number of
  polynomials of degree at most $d$ in $k$ variables with bitsizes bounded by
  $\tau$ is bounded by $(2^{\tau + 1})^{\binom{d + k}{k}} = 2^{\tau d^{O
  (k)}}$.
\end{remark}

Theorem \ref{13:the:saball} gives an upper bound on the radius of a ball meeting every connected
component of the realizations of every realizable sign condition on a family of polynomials with
integer coefficients. It is well known that any the intersections of any semi-algebraic set $S \subset \R^k$
with open balls of large enough radius are semi-algebraically homeomorphic. This is a consequence of 
the  local conic structure of
semi-algebraic sets \cite[Theorem 9.3.6]{BCR}.
The following theorem gives a single
exponential upper bound on the radius of a ball such that the \emph{homotopy type} (not necessarily the homeomorphism type) of a given semi-algebraic set defined by polynomials with coefficients in $\Z$,
is preserved upon intersection of with all balls of larger radii.

\begin{theorem}\cite{BV06}
\label{thm:ball}
Let $\mathcal{P}=\{P_1, \ldots, P_s \} \subset
  \mathbb{Z}[X_1, \ldots, X_k]$ and suppose that $P \in \mathcal{P}$ have
  degrees at most $d$, and the coefficients of $P \in \mathcal{P}$ have
  bitsizes at most $\tau$. 
Let $S \subset {\R}^k$ be a ${\mathcal P}$-semi-algebraic set.

There exists a constant $c > 0$, such that for any
$R_1 > R_2 > 2^{\tau d^{ck}}$ we have,
\begin{enumerate}
\item
$S \cap \B_k(0,R_1) \simeq S \cap \B_k(0,R_2)$, and
\item
$S \setminus \B_k(0,R_1) \simeq S \setminus \B_k(0,R_2)$
\end{enumerate}
(where $\simeq$ denotes semi-algebraic homotopy equivalence).
\end{theorem}

\begin{theorem}\cite{BV06}
\label{thm:contraction}
Let $S \subset {\R}^m$ be a ${\mathcal P}$-semi-algebraic set,
with ${\mathcal P} \subset \Z [X_1, \ldots , X_k]$ and $\mathbf{0} \in S$.
Let $\deg (P) <d$ for each $P \in {\mathcal P}$, and
the  bitsizes of the coefficients of $P \in {\mathcal P}$
be less than $\tau$.
Then, there exists a constant $c >0$ such that
for every $0 < r < 2^{- \tau d^{ck}}$ the set $S \cap \B_k(0,r)$ is contractible.
\end{theorem}

\begin{question}
Can one replace the relation of ``semi-algebraic homotopy equivalence''  in 
Theorem \ref{thm:ball}  by ``semi-algebraic homeomorphism'' while preserving the 
same bound? 
It follows from effective versions of Hartd's triviality theorem
\cite{Hardt} (see also \cite[Theorem 9.3.2.]{BCR})
that one can make this replacement, but with a doubly exponential bound on the radius --
i.e. a bound of the form $2^{\tau 2^{d^{ck}}}$ rather than $2^{\tau d^{ck}}$. A singly exponential bound
as in Theorem \ref{thm:ball} would be of much greater interest.
\end{question}

Another result along the lines of Theorem  
\ref{13:the:saball}
that has several 
for example in constructing certificates of positivity of polynomials 
in simplices using effective versions of Polya's theorem
(see \cite{Powers-Reznick-2001})
is the following theorem.

\begin{theorem} \cite{JP2011}
\label{thm:simplex}
Let $P \in  \mathbb{Z}[X_1, \ldots, X_k]$ 
with $\deg(P) \leq d$,  and the coefficients of $P$ have  
bitsizes bounded by $\tau$, and such that $P$ is strictly positive in the 
closed standard simplex  $\mathbf{\Delta}_k \subset \mathbb{R}^k$. Then, the 
minimum value that $P$ attains in $\mathbf{\Delta}_k$ is bounded from below by
\[
2^{-(\tau + 1) d^{k+1}} d^{-(k+1)d^{k}} \binom{d+1}{k+1}^{d^k(d-1)}.
\]
\end{theorem}

A related result is the following bound on the minimum value attained by an integer polynomial 
restricted to a compact connected component of a basic closed semi-algebraic subset of $\mathbb{R}^k$
defined by polynomials with integer coefficients in terms of the degrees and the bitsizes of the coefficients
of the polynomials involved. 

\begin{theorem} \cite{GPT2013}
\label{thm:opt}
Let $\mathcal{P} = \{P_1,\ldots.P_s\}  \subset  \mathbb{Z}[X_1, \ldots, X_k]$
with $\deg(P) \leq d$ for all $P \in \mathcal{P}$, and let the coefficients of $P$ have  
bitsizes bounded by $\tau$.  Let $Q \in \Z[X_1,\ldots,X_k]$, $\deg(Q) \leq d$, and
let the bit-sizes of the coefficients of $Q$ be also bounded by $\tau$.
Let $C$ be a compact connected component of the basic closed semi-algebraic set defined by
$P_1 = \cdots = P_\ell = 0, P_{\ell+1} \geq 0,\ldots, P_s \geq 0$. Then the minimum value 
attained by $Q$ over $C$ is a real algebraic number of degree at most $2^{n-1} d^n$,
and if it is not equal to $0$, then its absolute value is bounded from below by 
\[
(2^{4 - k/2}H d^k)^{-k 2^k d^k},
\]
where $H = \max(2^\tau, 2k + 2s)$.
\end{theorem}

\subsection{Singly exponential quantifier elimination algorithms}
The algorithm mentioned in Theorem \ref{13:the:samplealg}, which has singly exponential 
complexity and 
produces sample points in every semi-algebraically
connected component of the realization of every realizable sign condition of a family
of polynomials,  is an important ingredient for designing efficient algorithm
for quantifier elimination.  
More precisely, used in a parametrized form, it allows us to eliminate
one whole block of variables (quantified by the same quantifier) in one step, 
unlike in algorithms based on cylindrical algebraic decomposition,
where the elimination has to proceed one variable at a time {\em regardless
of the block structure of the quantifiers}. 
The singly exponential algorithm
for eliminating one block of variables at a time is formalized as the
\bfdef{Block Elimination Algorithm} \cite[Chapter 14]{BPRbook2} 
and does the following.
Given a finite family of polynomials $\mathcal{P} \subset 
\R[X_1,\ldots,X_k,Y_1,\ldots,Y_\ell]$, the Block Elimination Algorithm 
produces as output a family of polynomials 
$\BElim_X(\mathcal{P}) \subset \R[Y_1,\ldots,Y_\ell]$.
The family $\BElim_X(\mathcal{P})$ has the following important property
that justifies its name.
For each
semi-algebraically connected component, $C \subset \R^\ell$,  
of each realizable sign condition of $\BElim_X(\mathcal{P})$, 
the set of realizable sign conditions of 
$\mathcal{P}(y) \subset \R[X_1,\ldots,X_k]$ stay invariant as $y$ is allowed to
vary over $C$. The Block Elimination Algorithm also produces a set of
parametrized (by $y$) sample points which are guaranteed to meet each 
semi-algebraically connected component of
the set of realizable sign conditions of 
$\mathcal{P}(y) \subset \R[X_1,\ldots,X_k]$.
The complexity of
this algorithm is bounded by $s^{k+1} d^{O(\ell+k)}$, where as usual
$s = \card \; \mathcal{P}$ and $d$ is a bound on the degrees of the polynomials
in $\mathcal{P}$.

\subsubsection{Sign Determination Algorithm}
\label{sec:BKR}
The {\em Block Elimination Algorithm} is one important ingredient of the
critical point based quantifier elimination algorithm. The other important 
ingredient is a \bfdef{Sign Determination Algorithm}
that allows one to compute
the vector of signs of a family, $\mathcal{P}$, 
of $s$ polynomials in $\D[X]$ at the real 
roots of a fixed polynomial $Q \in \D[X]$, with complexity $s d^{O(1)}$,
where $d$ is a bound on the degrees of the polynomials in  
$\mathcal{P}$ and $Q$. This algorithm was first discovered by Ben-Or, Kozen and
Reif \cite{BKR} and extended  by Roy and Szpirglas \cite{RS} (see also
\cite{Perrucci09} for recent improvements). This algorithm has also been
generalized to the multi-variate case (where the zeros of $Q$ could be
positive dimensional), and this is described below in Section \ref{sec:eq}.

\subsubsection{Quantifier Elimination Algorithm}
The above ingredients (namely, the  {\em Block Elimination Algorithm}
and the {\em Sign Determination Algorithm}), along with numerous
technical detail which we omit in this survey, 
allows one to prove the following result.

\begin{theorem}\cite{BPRbook2}
\label{14:the:tqe}
 Let
$\mathcal{P}$ be a set of at most $s$ polynomials each of degree at
 most $d$ in $k+\ell$ variables with
coefficients in a real closed field $\R$,
and let $\Pi$ denote a partition
 of the list of variables
$(X_1,\ldots, X_k)$ into blocks,
${X}_{[1]},\ldots,{X}_{[\omega]}$, where the block $X_{[i]}$
has size $k_i, 1 \leq i \leq
\omega$.
Given $\Phi(Y)$,
 a  $(\mathcal{P},\Pi)$-formula, there exists an equivalent
quantifier free formula,
\[ \Psi(Y) =\bigvee_{i=1}^{I}\bigwedge_{j=1}^{J_i}(\bigvee_{n=1}^{N{i,j}}{\rm
sign}(P_{ijn}(Y))=
\sigma_{ijn}), \] where $P_{ijn}(Y)$ are polynomials in
the variables $Y$, $\sigma_{ijn} \in \{
0, 1 , -1 \},$
\[ I \leq s^{(k_\omega +1)\cdots(k_1 +1)(\ell+1)}d^{O(k_\omega) \cdots
O(k_1)O(\ell)}, \]
\[ J_i \leq s^{(k_\omega +1)\cdots(k_1 +1)}d^{O(k_\omega) \cdots O(k_1)} , \]
\[ N_{ij} \leq d^{O(k_\omega) \cdots O(k_1)} , \]
and the degrees of the polynomials $P_{ijn}(Y)$
are bounded by $d^{O(k_\omega) \cdots O(k_1)}$. Moreover,
there is
an algorithm to compute $\Psi(Y)$ with complexity
$$
s^{(k_\omega +1)\cdots(k_1 +1)(\ell  + 1)}d^{O(k_\omega) \cdots O(k_1)O(\ell)}
$$ in $\D$, denoting by
$\D$ the ring generated by the
coefficients of $\mathcal{P}$.

If $\D=\Z,$ and the bit-sizes
of the coefficients of the polynomials are bounded by
 $\tau$, then the bit-sizes of the integers appearing in the
intermediate computations and the output are bounded
by $\tau d^{O(k_\omega)\cdots O(k_1)O(\ell)}$.
\end{theorem}

\begin{remark}
\label{rem:critical}
The algorithmic results described in Section \ref{sec:critical}
are based on one common technique.  By taking
a well chosen infinitesimal perturbation, one can replace any bounded, real
(possibly singular) variety $V \subset \R^k$, by a non-singular variety
defined over an (non-archimedean) extension of $\R$. The 
projection map on some co-cordinate (say $X_1$) restricted to
this new variety has only non-degenerate critical points, which moreover are
defined by a zero-dimensional system of equations which is nicely behaved
(is automatically a Gr\"{o}bner basis). The limits of these critical points
belong to the given variety $V$ and moreover they meet every semi-algebraically
connected component of $V$. This technique (which is rather special to
real algebraic geometry as opposed to complex geometry) 
has several advantages from the point of view
of algorithmic complexity. The first advantage is that it is not necessary
to choose any generic co-ordinate system or direction to project on. Secondly,
the method does not care about how singular the given variety $V$ is 
or even its dimension. Moreover, it is possible to relate the topology 
(up to semi-algebraic homotopy equivalence) of $V$ with the infinitesimal
``tube'' around it which is bounded by the perturbed hyper-surface (say $V'$).
This reduces most algorithmic problems of computing topological invariants
of $V$, to that of the well-behaved hyper-surface $V'$. Since the degree
of the polynomial defining $V'$ is at most twice that of the one defining $V$,
and the computations take place in the original ring adjoined with at most
a constant many (i.e. their number is independent of the input parameters
$s,d$ and $k$)  infinitesimals, the complexity is well controlled. 
The main disadvantage of the approach (which could be a drawback from the
point of view of practical implementation) 
is that computations with even a constant many
infinitesimals are quite expensive (even though they do not affect
the asymptotic complexity bounds). Also, the process of taking algebraic
limits at the end can  be quite cumbersome. Nevertheless, this perturbation
approach remains the only one which gives deterministic algorithms with the
best known worst case complexity estimates.
\end{remark}

\subsection{Intrinsic complexity and complex algebraic techniques}
The model for studying complexity of algorithms in this survey is that
the size of the input is measured in terms of the number of coefficients
needed to specify the input polynomials in the dense representation. 
Since this number is determined by the following parameters:
\begin{enumerate}
\item
the number of variables, $k$;
\item
the number of polynomials, $s$;
\item
the degrees of the polynomials, $d$;
\end{enumerate}
it makes sense to state the complexity estimates in terms of $s,d$ and $k$.

There is another body of work (see for example \cite{BGHM97,BGHM01,MS03,MS04,
JPS09}) 
in which the goal is to obtain algorithms for computing sample points
on each semi-algebraically connected component of a given 
real algebraic variety $V \subset \R^k$,  
whose complexity is bounded by a {\em polynomial} function of 
some {\em intrinsic}
invariant of the variety $V$ or in some cases the  length of  
{\em straight line programs} encoding  the input polynomials. 
In this approach, the real variety $V$ is considered as the real part of
the complex variety $V_\C \subset \C^k$ (where $\C$ is the algebraic closure
of $\R$), and the intrinsic invariant,
$\delta(V) = \delta(V_\C)$ 
depends only on the geometry of the complex variety $V_\C$,
and not on the
particular presentation of it by the given input polynomials. 
If $d$ is a bound on the degrees of the polynomials defining $V$, then
$\delta(V)$ is bounded by $O(d)^k$ and could be as large as
$d^k$ in the  worst case. However, $\delta(V)$ could be smaller in
special cases.

Since these algorithms aim at complexity in terms of some geometric invariant
of the variety itself, the infinitesimal perturbation techniques described
in the previous sections is not available, 
since such a perturbation will not in general
preserve this invariant. Hence, one needs to work directly with the given 
variety. 
For example, one needs to prove that under certain assumptions on
the variety, the critical points of a generic projection 
(also called the polar variety)
is non-singular (see \cite{BGHMS10}).  
The theory of {\em geometric resolutions} (see \cite{BGHM01}) play
an important role in these algorithms.

One feature of the algorithms that follow from these techniques
is that it is necessary  to choose generic co-ordinates which cannot be done
deterministically within the claimed complexity bounds. As such one obtains
{\em probabilistic} (as opposed to deterministic) algorithms, meaning that
these algorithms always run within the stated complexity time bounds, 
but is guaranteed to give correct results only with high probability. 

\subsection{Variants of quantifier elimination and applications}
In certain applications (most notably in the theory of constraint
databases) one needs to perform quantifier elimination in a more
generalized setting than that discussed 
above.
For instance,
it is sometimes necessary to eliminate quantifiers not just from one
formula, but a whole sequence of formulas described in some finite
terms, where the number of free variables is allowed to grow in the
sequence. Clearly, the quantifier elimination algorithms described previously
is not sufficient for this purpose since their 
complexity depends on the number of free variables.

We describe below 
a variant of the quantifier elimination problem
which was introduced in \cite{B99b} motivated by
a problem in constraint databases.

\subsubsection{The Uniform Quantifier Elimination Problem}
\begin{definition}
\label{def:uniformqe}
We call a sequence, 
\[
\{\phi_n(T_1,\ldots,T_\ell,Y_1,\ldots,Y_n) \;|\; n > 0 \}
\]
of first-order formulas $\phi_n$ in the language of ordered fields,
to be a \bfdef{uniform sequence}  if each $\phi_n$ has the form,
\[ 
\phi_n(T_1,\ldots,T_\ell,Y_1,\ldots,Y_n) = 
\]
\[
Q^1_{1 \leq k_1 \leq n}\ldots
   Q^{\omega}_{1 \leq k_{\omega} \leq n} \phi(T_1,\ldots,T_\ell,Y_{k_1},\ldots,
Y_{k_{\omega}}) ,
\]
where $Q^i \in \{\vee,\wedge\}, 1 \leq i \leq \omega $ and 
$\phi$ is some  fixed $(\ell+\omega)$-ary quantifier-free 
first-order formula.

Thus for every $n,$  $\phi_n$ is a first order formula with $\ell+n$
free variables. We will refer to the variables $T_1,\ldots,T_\ell$ as
{\em parameters}.

Given a uniform sequence of formulas
$
\Phi = \{ \phi_n \;|\; n > 0 \},
$
where 
\[ 
\phi_n(T_1,\ldots,T_\ell,Y_1,\ldots,Y_n) =  
\]
\[
Q^1_{1 \leq k_1 \leq n}\ldots
   Q^{\omega}_{1 \leq k_{\omega} \leq n} \phi(T_1,\ldots,T_\ell,Y_{k_1},\ldots,
Y_{k_{\omega}}) ,
\]
we define the {\em size} of $\Phi$ 
to be the
length of the formula $\phi$. 
\end{definition}

\begin{example}
\label{eg:uniformqe}
Consider the uniform sequence of formulas
\[ 
\phi_n(T_1,Y_1,\ldots,Y_n) = \bigwedge_{1 \leq k_1 \leq n} (Y_{k_1} - T_1 = 0),
\;\; n > 0 .
\]

Consider the  sequence of quantified formulas, 
$(\exists T_1) \phi_n(T_1,Y_1,\ldots,Y_n) .$ 
In this example, it is easily seen that letting
\[
\Psi_n =  \bigwedge_{1 \leq k_1 \leq n} \bigwedge_{1 \leq k_2 \leq n}
(Y_{k_1} - Y_{k_2} = 0), 
\]
we get a uniform sequence of quantifier-free formulas satisfying,
\[
\Psi_n(Y_1,\ldots,Y_n) \Leftrightarrow (\exists T_1) \phi_n(T_1,Y_1,\ldots,Y_n)
\]
for
every $n > 0$.
\end{example}

The \bfdef{uniform quantifier elimination problem}  is
to eliminate quantifiers from a uniform sequence
of formulas and obtain another {\em uniform} sequence of quantifier
free formulas. 

The following is proved in \cite{B99b}.

\begin{theorem}(Uniform Quantifier Elimination)
\label{the:uniformqe}
Let, 
\[
\Phi = \{\phi_n(T_1,\ldots,T_\ell,Y_1,\ldots,Y_n) \;|\; n > 0 \}
\]
be a uniform sequence of formulas with  parameters $T_1,\ldots,T_\ell,$
where 
\[ 
\phi_n(T_1,\ldots,T_\ell,Y_1,\ldots,Y_n) =   
\]
\[
Q^1_{1 \leq k_1 \leq n}\ldots
   Q^{\omega}_{1 \leq k_{\omega} \leq n} \phi(T_1,\ldots,T_\ell,Y_{k_1},\ldots,
Y_{k_{\omega}}).
\]

Let the number of different $(\ell+\omega)$-variate 
polynomials appearing in $\phi$ be
$s$ and let their degrees be bounded by $d$.

Let $R_1,\ldots, R_m \in \{\exists,\forall \},$
$R_i \neq R_{i+1},$ and let 
$T^{[1]},\ldots,T^{[m]}$ be a partition of the variables,
$T_1,\ldots,T_\ell$ into $m$ blocks of size $\ell_1,\ldots,\ell_m,$
where $\sum_{1 \leq i \leq m} \ell_i = \ell$.

Then, there exists an algorithm that outputs a quantifier-free
first order formula, 
$
\psi(Y_{k_1},\ldots,Y_{k_{\omega'}}),
$
along with $Q^i \in \{\bigvee,\bigwedge \}, 1 \leq i \leq \omega',$
such that for every $n >0$
\[ 
\psi_n(Y_1,\ldots,Y_n) = 
Q^1_{1 \leq k_1 \leq n}\ldots
   Q^{\omega'}_{1 \leq k_{\omega'} \leq n} \psi(Y_{k_1},\ldots,
Y_{k_{\omega'}}) 
\] 
\[
\Leftrightarrow
(R_1 T^{[1]})\ldots(R_m T^{[m]})
\phi_n(Y_1,\ldots,Y_n,T_1,\ldots,T_\ell).
\] 

The complexity of the algorithm is bounded by
\[
s^{\prod_i (\ell_i +1)} d^{\omega\prod_i O(\ell_i^2)},
\]
and the size of the formula $\psi$ is bounded by 
\[
s^{\prod_i (\ell_i +1)} d^{\omega\prod_i O(\ell_i^2)}\mathrm{size}(\phi).
\]
\end{theorem}

\begin{remark}
\label{rem:inexpressibility}
In \cite{B99b}
Theorem \ref{the:uniformqe} is used to prove the equivalence of
two different semantics and in the theory of constraint databases. However,
it also has applications in logic. For example, in the same paper it is
used to prove that semi-algebraic connectivity is not expressible by a 
first-order formula (see \cite{B99b} for a precise definition of first-order
expressibility). 
Note that the inexpressibility result was also proved by more abstract
model theoretic methods in \cite{Libkin98}.

The technique used in the proof of Theorem \ref{the:uniformqe} 
is also used in \cite{B99b}
to give an algorithm for ordinary quantifier elimination 
whose complexity depends on the size of the input formula, 
and which has better complexity than the
algorithm in Theorem \ref{14:the:tqe} 
in case the input formula has a small size.
This algorithm is called \bfdef{Local Quantifier Elimination Algorithm}
in \cite{BPRbook2}.
\end{remark}

\section{Computing topological invariants of semi-algebraic sets}
\label{sec:top}
As remarked above (see Remark \ref{rem:inexpressibility}), an effective 
algorithm for deciding connectivity of semi-algebraic sets does not
automatically follow from the Tarski-Seidenberg principle. However,
one can decide questions about connectivity (as well as compute other 
topological invariants such as the Betti numbers) using effective
triangulation of semi-algebraic sets via Cylindrical
Algebraic Decomposition. However, such an algorithm will necessarily have 
doubly exponential complexity.

Most of the recent work in algorithmic  semi-algebraic
geometry has focused on obtaining {\em singly exponential time}
algorithms -- that is algorithms with complexity of the
order of 
$(s d)^{k^{O(1)}}$ 
rather than $(s d)^{2^k}$. 
An important motivating reason behind the search for such algorithms,
is the following theorem  due to Gabrielov and Vorobjov \cite{GV07}
(see also \cite{GaV}) 
(see  \cite{OP,T,Milnor2,Basu1}, as well as the 
survey article \cite{BPR10}, for work leading up to this result)
which gives singly exponential upper bound on the topological
complexity of semi-algebraic sets measured by the sum of
their Betti numbers.

\begin{theorem} \cite{GV07}
\label{the:GV}
For a ${\mathcal P}$-semi-algebraic set $S \subset \R^k$,
the sum of the Betti numbers of $S$ 
is bounded by $(O(skd))^k$,
where 
$s = \card \;{\mathcal P}$, and $d = \max_{P \in {\mathcal P}} \deg(P)$.
\end{theorem}

For the special case of ${\mathcal P}$-closed semi-algebraic sets the
following slightly better bound was known before \cite{Basu1}
(and this bound is used in an essential way in the proof of
Theorem \ref{the:GV}).
Using the same notation as in Theorem \ref{the:GV} above we have

\begin{theorem} \cite{Basu1}
\label{the:B99}
For a ${\mathcal P}$-closed semi-algebraic set $S \subset \R^k$,
the sum of the Betti numbers of $S$ is bounded by $(O(sd))^k$.
\end{theorem}

\begin{remark}
\label{rem:lowerbound}
These bounds are asymptotically tight, as can be already seen from the 
example where each $P \in {\mathcal P}$ is a product of 
$d$ generic polynomials of degree one. The number of semi-algebraically connected
components of the ${\mathcal P}$-semi-algebraic
set defined as the subset of $\R^k$
where all polynomials in ${\mathcal P}$ are non-zero is clearly
bounded from below by $(C sd)^k$ for some constant $C$.
\end{remark}

\subsection{Roadmaps}
\label{sec:roadmap}
Theorem \ref{13:the:samplealg} gives a singly exponential time
algorithm for testing if a given semi-algebraic set is empty or not.
However, it gives no way of testing if any two sample points computed
by it belong to the same semi-algebraically connected component of the given
semi-algebraic set, even though the set of sample points is guaranteed to meet
each such semi-algebraically connected component. 
In order to obtain connectivity information in singly exponential time
a more sophisticated construction is required -- namely that of a {\em roadmap}
of a semi-algebraic set, which is an one dimensional semi-algebraic subset
of the given semi-algebraic set which is non-empty and semi-algebraically connected inside
each semi-algebraically connected component of the given set. Roadmaps
were first introduced by Canny \cite{Canny93a}, but similar constructions
were considered as well by Grigoriev and Vorobjov \cite{GV92} and
Gournay and Risler \cite{GR92}. 
Our exposition below follows that in \cite{BPR99,BPRbook2} where the
most efficient algorithm for computing roadmaps is given.
The notions of pseudo-critical points and values defined above play a
critical role in the design of efficient algorithms for computing 
roadmaps of semi-algebraic sets.

We first define a \bfdef{roadmap  of a semi-algebraic set}. 
We use the following notation.
We  denote by $\pi_{1\ldots j}$ the projection,
$x \mapsto (x_1,\ldots,x_j).$
Given a set $S \subset \R^k$ and $y \in \R^j$, we denote by 
$S_y=S \cap \pi_{1 \ldots j}^{-1}(y)$.

\begin{definition}[Roadmap of a semi-algebraic set]
\label{15:def:roadmap}
Let $S \subset \R^k$ be a semi-algebraic set.
A {\em roadmap}
\index{Roadmap} for $S$
 is a semi-algebraic set $M$
of dimension at most one contained in $S$
which satisfies the following roadmap
conditions:
\begin{itemize}
\item ${\rm RM}_1$ For every semi-algebraically
connected component $D$ of $S$,
$D \cap M$ is non-empty and semi-algebraically connected.
\item ${\rm RM}_2$ For every $x \in {\R}$ and
for every semi-algebraically connected component $D'$
of $S_x$, $D'\cap  M \neq \emptyset$.
\end{itemize}
\end{definition}

We describe the construction of a roadmap
$\RM(\ZZ(Q,\R^k),{\mathcal N})$ for a bounded algebraic set $\ZZ(Q,\R^k)$
 which contains a finite set of
points ${\mathcal N}$ of $\ZZ(Q,\R^k)$. A precise description of how the
construction can be performed algorithmically
can be found in \cite{BPRbook2}.
We should emphasize here that $\RM(\ZZ(Q,\R^k),{\mathcal N})$ denotes the
semi-algebraic set output by the specific algorithm described below which
satisfies the properties stated in Definition \ref{15:def:roadmap}
(cf. Proposition \ref{15:prop:rm}).

Also, in order to understand the roadmap algorithm it is easier to first
concentrate on the case of a bounded and non-singular real algebraic
set in $\R^k$ (see Figure \ref{fig:torus2} below). In this case several
definitions get simplified. For example, the pseudo-critical values
defined below are in this case ordinary critical values of the projection map
on the first co-ordinate. However, one should keep in mind that even 
if one starts with a bounded non-singular algebraic set, the input to
the recursive calls corresponding to the critical sections (see below)
are necessarily singular and thus it is not possible to treat the 
non-singular case independently.
\begin{center}
\begin{figure}
\includegraphics[height=7cm]{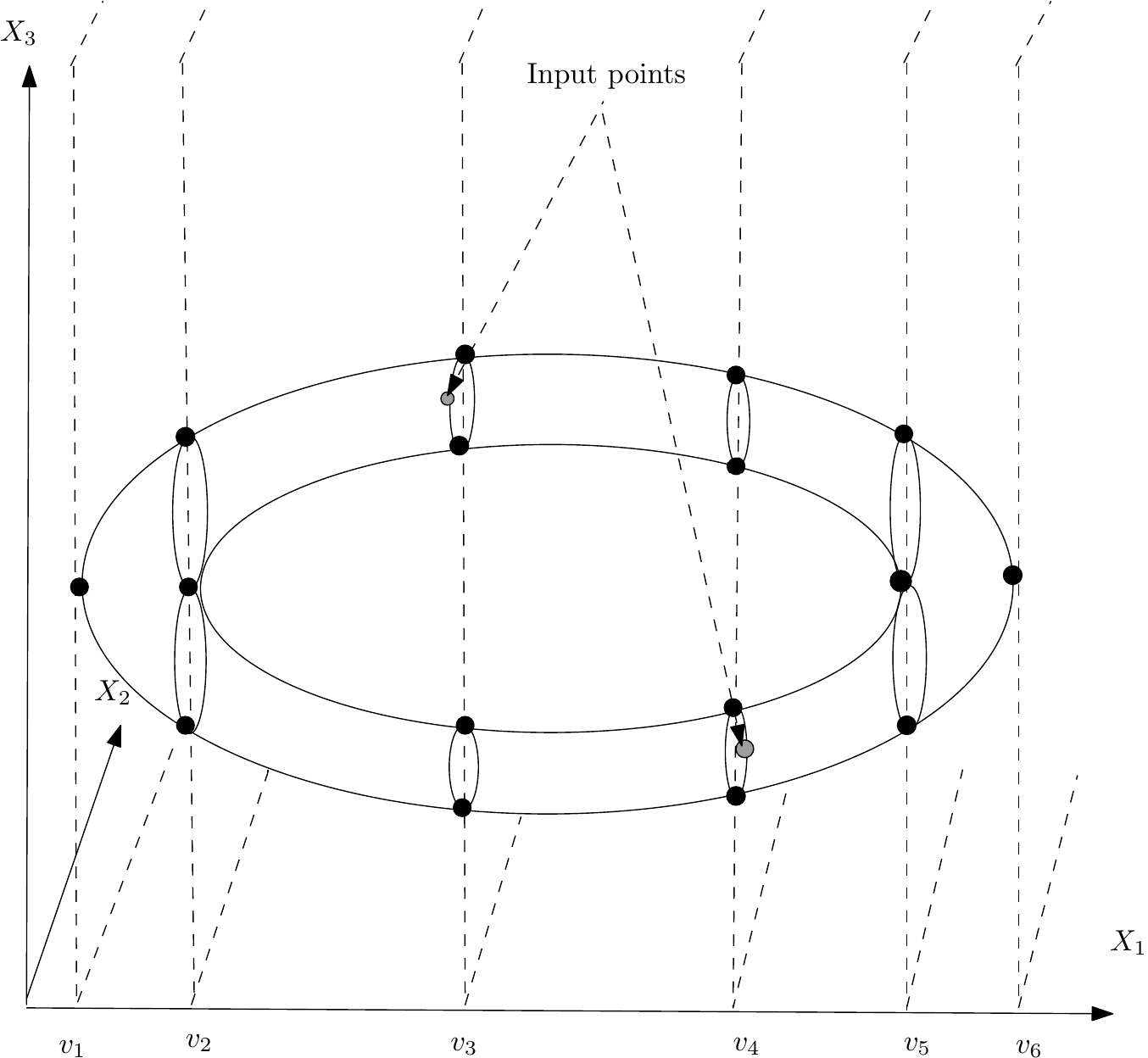}
\caption{Roadmap of the torus in $\R^3$.}
\label{fig:torus2}
\end{figure}
\end{center}
A key ingredient of the roadmap is  the construction of 
pseudo-critical points and values defined above.
The construction of the roadmap of an algebraic set
containing a finite number 
of input points ${\mathcal N}$ of this algebraic set
is as follows.
We first construct $X_2$-pseudo-critical points on $\ZZ(Q,\R^k)$
in a parametric way along  the $X_1$-axis by following continuously, 
as $x$ varies on the $X_1$-axis,
the $X_2$-pseudo-critical points on $\ZZ(Q,\R^k)_{x}$. This results in curve segments
and their endpoints on $\ZZ(Q,\R^k).$ The curve segments
 are continuous semi-algebraic curves parametrized by open intervals on the 
$X_1$-axis and their endpoints are points of $\ZZ(Q,\R^k)$
above the corresponding endpoints of the open intervals.
Since these curves and their endpoints include
for every $x\in\R$ the $X_2$-pseudo-critical points of
$\ZZ(Q,\R^k)_{x}$, they meet every semi-algebraically connected component of
$\ZZ(Q,\R^k)_{x}$.  Thus, the set of curve
segments and their endpoints already satisfies
${\rm RM}_2.$ However, it is clear that this set might not
be semi-algebraically connected in a semi-algebraically
connected component and so ${\rm RM}_1$ might not be satisfied.
We add additional curve segments to ensure connectedness by recursing 
in certain  distinguished hyper-planes defined by
$X_1=z$ for distinguished values $z$.

The set of {\em distinguished values}
is the union of the $X_1$-pseudo-critical values, the first
coordinates of the input points ${\mathcal N}$, and the
first coordinates of the endpoints of the curve
segments. A {\em distinguished hyper-plane}
is an hyper-plane defined by $X_1=v$, where
$v$ is a distinguished value. The input points, the endpoints of
the curve segments, and the intersections of the curve
segments with the distinguished hyper-planes define 
the set of {\em distinguished points}.

Let
the distinguished values be
$v_1<\ldots <v_\ell.$
Note that amongst these  are the $X_1$-pseudo-critical values. Above each
interval $(v_i, v_{i+1})$ we have  constructed a collection
of curve segments ${\mathcal C}_i$ meeting every
semi-algebraically connected component of
$\ZZ(Q,\R^k)_v$ for every $v \in (v_i, v_{i+1})$. Above
each distinguished value $v_i$  we have 
a set of distinguished points ${\mathcal N}_i$. 
Each curve segment in ${\mathcal C}_i$
has an endpoint in ${\mathcal N}_i$ and another in
${\mathcal N}_{i+1}$.  Moreover, the union of the ${\mathcal N}_i$
contains ${\mathcal N}$.

We then repeat this construction in each 
distinguished hyper-plane $H_i$ defined by $X_1=v_i$
with input $Q(v_i,X_2,\ldots,X_k)$ and the
distinguished points in ${\mathcal N}_i$. 
Thus, we construct distinguished values
$v_{i,1},\ldots, v_{i,\ell(i)}$ of 
$\ZZ(Q(v_i,X_2,\ldots,X_k),\R^{k-1})$ 
(with the role of $X_1$ being now played by $X_2$) and 
the process is iterated until
for 
$I=(i_1,\ldots,i_{k-2}), 1 \leq i_1 \leq\ell,
\ldots, 1 \leq i_{k-2} \leq \ell(i_1,\ldots,i_{k-3}),$ we have
distinguished values
$v_{I,1}< \ldots < v_{I, \ell(I)}$ along the
$X_{k-1}$ axis with corresponding sets of curve
segments and sets of distinguished points with the
required incidences between them.

The following theorem  is proved in \cite{BPR99}
(see also \cite{BPRbook2}).
\begin{proposition}
\label{15:prop:rm}
The semi-algebraic set
$\RM(\ZZ(Q,\R^k),{\mathcal N})$ obtained by this construction is a
roadmap for $\ZZ(Q,\R^k)$ containing ${\mathcal N}$.
\end{proposition}

Note that if $x \in \ZZ(Q,\R^k)$,  then $\RM(\ZZ(Q,\R^k),\{x\})$ 
contains a path,
$\gamma(x)$,  connecting a distinguished point $p$ 
of   $\RM(\ZZ(Q,\R^k))$ to $x$.

\subsubsection{Roadmaps of general semi-algebraic sets}
\label{subsecsec:generalroadmap}
Using the same ideas as above and some additional techniques for
controlling the combinatorial complexity of the algorithm 
it is possible to extend 
the roadmap algorithm to the case of semi-algebraic sets. The following
theorem appears in \cite{BPR99,BPRbook2}.

\begin{theorem}\cite{BPR99,BPRbook2}
\label{16:the:saconnecting}
Let $Q \in
{\R}[X_1,\ldots,X_k]$
 with $\ZZ(Q,\R^k)$  of dimension
$k'$ and let ${\mathcal
P}\subset {\R}[X_1,\ldots,X_k]$
be a set of at most $s$
polynomials
for
which the degrees of the polynomials
 in ${\mathcal P}$ and $Q$ are bounded by
$d.$
Let $S$ be a ${\mathcal P}$-semi-algebraic subset of $\ZZ(Q,\R^k)$.
There is an algorithm which computes a roadmap $\RM(S)$ for $S$
with complexity $s^{k'+1} d^{O(k^2)}$
 in the ring ${\D}$
generated by the
coefficients of $Q$ and the elements
of ${\mathcal P}$.
If $\D=\Z,$ and the
bit-sizes
of the coefficients of the polynomials are bounded by
 $\tau$,
then the bit-sizes of the integers appearing in the
intermediate
computations and the output are bounded
by $\tau  d^{O(k^2)}$.
\end{theorem}

Theorem \ref{16:the:saconnecting} immediately implies that
there is an algorithm whose output is exactly one point
in every semi-algebraically connected component of
$S$ and  whose complexity in the ring generated by the 
coefficients of $Q$ and ${\mathcal P}$ is bounded by
$s^{k'+1} d^{O(k^2)}$.
In particular, this algorithm counts
the number semi-algebraically connected component of
$S$ within the same time bound.

\subsubsection{Recent  developments}
\label{sec:baby-giant}
Very recently Schost and Safey el Din \cite{Mohab-Schost2010} have
given a {\em probabilistic} algorithm for computing the roadmap of a 
smooth, bounded real algebraic hyper-surface in $\R^k$ 
defined by a polynomial of degree $d$, whose complexity is bounded
by $d^{O(k^{3/2})}$. 
Complex algebraic techniques related to the geometry of polar varieties
play an important role in this algorithm.
More recently, a {\em deterministic} 
algorithm for computing 
roadmaps of {\em arbitrary} 
real algebraic sets with the same complexity bound,
has also been obtained
\cite{BRMS10}. 

This algorithm is based on techniques coming from semi-algebraic
geometry and can be seen as a direct generalization of 
Proposition \ref{15:prop:rm} above. The main new idea is to consider
the critical points of projection maps onto a co-ordinate subspace of
dimension bigger than $1$ (in fact, of dimension $\sqrt{k}$). As a result
the dimensions in the recursive calls to the algorithm decreases by
$\sqrt{k}$ at each step of the recursion
(compared to the case of the ordinary roadmap algorithms where 
it decreases by $1$ in each step). This results in the improved complexity.
One also needs to prove suitable generalizations of the results 
guaranteeing the connectivity of the roadmap (see \cite[Chapter 15]{BPRbook2})
in this more general situation.
 
The recursive schemes used in the algorithms in \cite{Mohab-Schost2010} and \cite{BRMS10} have a 
commmon defect in that they are unbalanced in the following sense. 
The dimension of the fibers in which recursive calls is equal to $k -\sqrt{k}$
(in the classical case this dimension is $k-1$), which is much larger than
the dimension of the ``polar variety'' which is $\sqrt{k}$ (this dimension is equal to $1$ in the classical case).
While being less unbalanced than in the classical case (which accounts for the improvement in the complexity), there is scope for further improvement if these two dimensions can be made roughly equal.
There are certain formidable technical obstructions to be overcome to achieve this. 

This was done in  
\cite{BR14} 
where an algorithm based on a balanced  (divide-and-conquer) scheme is given for computing 
a roadmap of an algebraic set. The following theorem is proved.

\begin{theorem}
\cite{BR14}
\label{thm:divide-and-conquer}
 Let $\R$ be a real closed field and $\D \subset \R$ an
ordered domain. 
There exists an algorithm that takes as input:
\begin{enumerate}
\item a polynomial $P \in \D [ X_{1} , \ldots ,X_{k} ]$, with $\deg ( P
) \leq d$;
\item a finite set, $A$, of real univariate representations
whose associated set of points, $\mathcal{A} = \{ p_{1} , \ldots ,p_{m}
\}$, is contained in $V= \ZZ ( P, \R^{k} )$, and such that
the degree of the real univariate representation representing $p_{i}$ is
bounded by $D_{i}$ for $1 \leq i \leq m$;
\end{enumerate}
and computes a roadmap of $V$ containing $\mathcal{A}$. The complexity of
the algorithm is bounded by 
\[\left( 1+ \sum_{i=1}^{m} D_{i}^{O ( \log^{2}
( k ) )} \right) ( k^{\log ( k )} d )^{O ( k\log^{2} ( k ) )}.
\] 
The
size of the output is bounded by $( \mathrm{card} ( \mathcal{A} ) + 1 ) (
k^{\log ( k )} d )^{O ( k\log ( k ) )}$, while the degrees of the
polynomials appearing in the descriptions of the curve segments and points
in the output are bounded by 
\[
( \max_{1 \leq i \leq m}D_{i} )^{O ( \log
( k ) )} ( k^{\log ( k )} d )^{O ( k\log ( k ) )}.
\]
\end{theorem}

A probabilistic algorithm based on a similar divide-and-conquer strategy which works for
smooth, bounded algebraic sets, and with a slightly better complexity is given in \cite{Mohab-Schost2014}. 
As in \cite{Mohab-Schost2010}, complex algebraic (as opposed to semi-algebraic) techniques
play an important role in this algorithm.

\begin{question}
Is it possible to extend the  algorithm in Theorem \ref{thm:divide-and-conquer} to general \emph{semi-algebraic} sets and thus obtain an algorithm for computing roadmaps of such sets having better complexity
than the one given in Theorem \ref{16:the:saconnecting}.
\end{question}

\subsubsection{Parametrized paths}
One important idea  in the algorithm for computing the
first Betti number of semi-algebraic sets, is   
the construction of certain semi-algebraic sets called 
{\em parametrized paths}.
Under a certain hypothesis, these sets  are semi-algebraically contractible.
Moreover, there exists an algorithm for computing 
a covering of a given basic semi-algebraic set, $S \subset \R^k$,
by a singly exponential number of parametrized paths. 

We are given  a polynomial $Q \in \R[X_1,\ldots,X_k]$
such that
$\ZZ(Q,\R^k)$ is bounded and
a finite set of polynomials ${\mathcal P} \subset \R[X_1,\ldots,X_k]$.

The main technical construction underlying the algorithm for
computing the first Betti number in \cite{BPRbettione}, is
to obtain a covering of a given 
$\mathcal P$-closed semi-algebraic set contained in $\ZZ(Q,\R^k)$
by a family of semi-algebraically contractible subsets.
This construction
is based on a parametrized version of the connecting algorithm:
we compute a family of polynomials such that for each realizable
sign condition $\sigma$ on this family, the description of the connecting
paths of different points in the realization, $\RR(\sigma,\ZZ(Q,\R^k)),$ are
uniform. 
We first define parametrized paths.
A parametrized path is a
semi-algebraic set which is a union of semi-algebraic
paths having a special property called  
the {\em divergence property} in \cite{BPRbettione}.

More precisely,
\begin{definition}[Parametrized paths]
\label{def:parametrizedpath}
A \bfdef{parametrized path}
$\gamma$ is
 a continuous semi-algebraic mapping  from 
$V \subset \R^{k+1}
\rightarrow \R^k,$
such that,
denoting by 
$U=\pi_{1\ldots k}(V)\subset \R^k$,
there exists a semi-algebraic continuous function
$\ell: U \rightarrow [0,+\infty),$ 
and there exists  a point $a$ in $\R^k$, such that
\begin{enumerate}
\item $V =\{(x,t) \mid x \in U, 0 \le t \le \ell(x)\},$
\item $\forall \; x \in U, \; \gamma(x,0)=a$,
\item $\forall \; x \in U, \; \gamma(x,\ell(x))=x$,
\item 
$$
\displaylines{
\forall \; x \in U, \forall \; y \in U, \forall \; s \in [0,\ell(x)], \forall \; t \in [0,\ell(y)]\cr
 \left(\gamma(x,s)=\gamma(y,t) \Rightarrow s=t \right),
}
$$
\item  
$$
\displaylines{
\forall \; x \in U, \forall \; y \in U, \forall \; s \in [0,\min(\ell(x),\ell(y))]\cr
\left(\gamma(x,s)=\gamma(y,s) \Rightarrow \forall \; t \le s \; \gamma(x,t)=\gamma(y,t) \right).
}
$$
\end{enumerate}
\end{definition}

\begin{center}
\begin{figure}
\includegraphics[height=7cm]{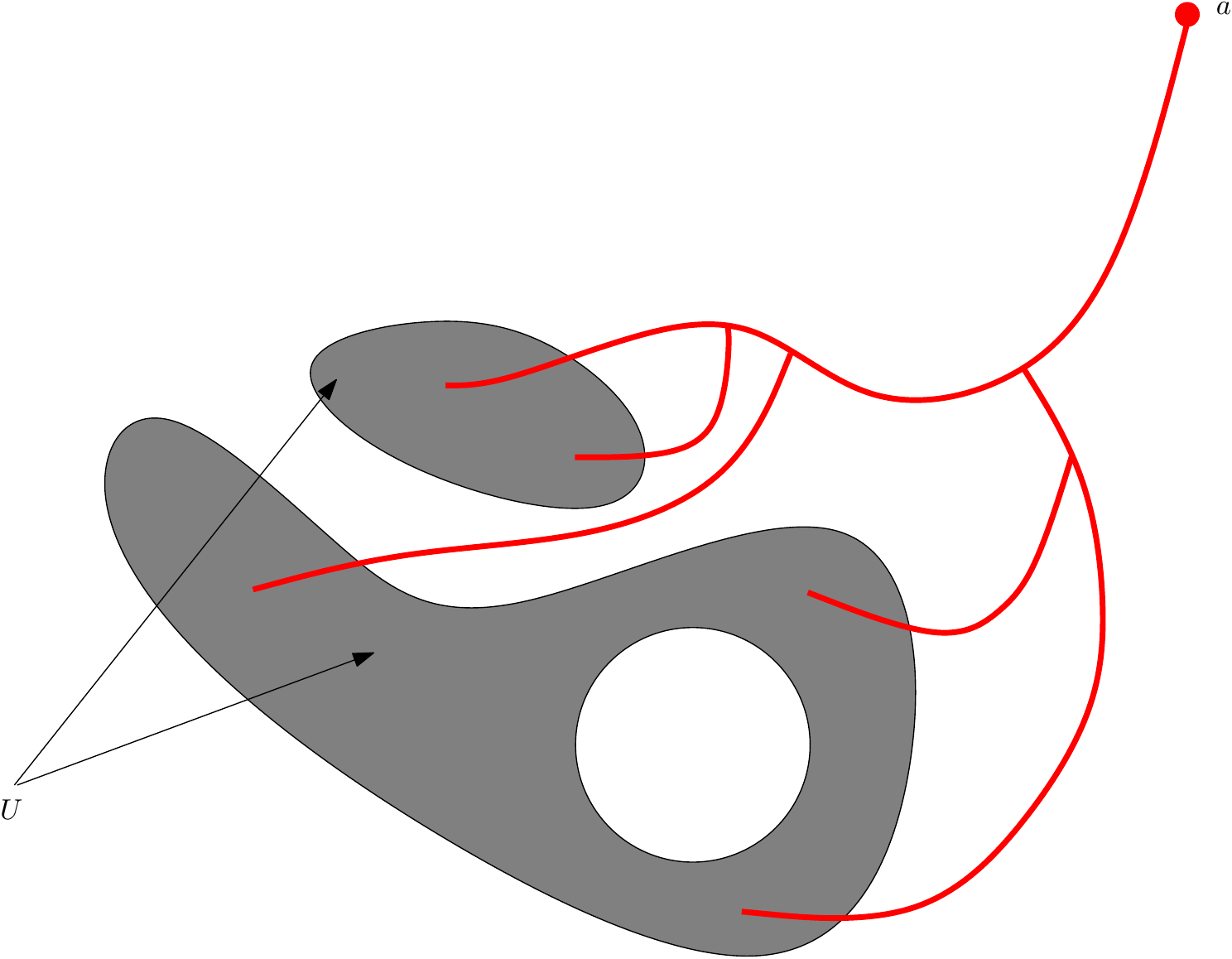}
\caption{A parametrized path}
\end{figure}
\end{center}

Given a parametrized path, $\gamma: V \rightarrow \R^k$, we will
refer to $U = \pi_{1\ldots k}(V)$ as its {\em base}.
Also,  any semi-algebraic subset 
$U' \subset U$ of the base of such a parametrized
path,
defines in a natural way  the restriction of 
$\gamma$ to the base $U'$, which is
another parametrized path,
obtained by restricting $\gamma$ to the set $V' \subset V$,
defined by 
$V' = \{(x,t) \mid x \in U', 0 \leq t \leq \ell(x) \}$.

The following proposition which appears in  \cite{BPRbettione} describes a crucial
property of parametrized paths, which makes them useful in algorithms
for computing Betti numbers of semi-algebraic sets.

\begin{proposition}\cite{BPRbettione}
\label{contractible}
Let $\gamma: V \rightarrow R^k$ be a parametrized path such that
$U = \pi_{1\ldots k}(V)$ is closed and bounded.
Then, the image of $\gamma$ is semi-algebraically contractible.
\end{proposition}

For every point $x$ of $\ZZ(Q,\R^k)$, denote by 
 $\sigma(x)$ the sign condition
on ${\mathcal P}$ at $x$. Let
$
 {\RR}(\overline\sigma(x),\ZZ(Q,\R^k)) =
\{x\in\ZZ(Q,\R^k)\mid\bigwedge_{P \in {\mathcal P}}\;\; \sign(P(x)) \in
\overline\sigma(x)(P)\},
$
where $\overline\sigma$ is the relaxation of $\sigma$
 defined by
$$
\left\{\begin{array}{ccc} \overline\sigma =
\{0\} &\mbox{ if }&\sigma =0,\\
\overline\sigma = \{0,1\} &\mbox{ if
}&\sigma =1,\\
\overline\sigma = \{0,-1\}& \mbox{ if }&\sigma
=-1.
\end{array}\right.
$$
We say that
$\overline \sigma(x)$ is the \bfdef{weak sign condition} defined by
 $x$ on ${\mathcal P}$. We denote by ${\mathcal P}(x)$
the union of $\{Q\}$ and the set of
polynomials in ${\mathcal P}$ vanishing at $x$.

The following theorem appears in \cite{BPRbettione}.

\begin{theorem}
\label{the:alg:parametrized}
There exists an algorithm that takes as input
a finite set of  polynomials ${\mathcal P} \subset \R[X_1,\ldots,X_k]$,
and produces as output,
\begin{itemize}
\item
a finite set of polynomials ${\mathcal A} \subset \R[X_1,\ldots,X_k]$,
\item 
a finite set $\Theta$ of quantifier free formulas, 
with atoms of the form $P = 0, P > 0, P< 0, \;P \in {\mathcal A}$,
such that for every semi-algebraically connected component $S$ of
the realization of
every weak sign condition on ${\mathcal P}$ on $\ZZ(Q, \R^k)$,
there exists a subset $\Theta(S) \subset \Theta$
 such that
$
\displaystyle{
S=\bigcup_{\theta \in \Theta(S)} \RR(\theta,\ZZ(Q,\R^k)),
}
$
\item for every $\theta \in \Theta$,
a parametrized path
$$
\displaylines{
\gamma_\theta : V_\theta \rightarrow \R^k,
}
$$
with base $U_\theta =  \RR(\theta,\ZZ(Q,\R^k))$, 
such that for each $y \in \RR(\theta,\ZZ(Q,\R^k)),$ 
the image of $\gamma_\theta(y,\cdot)$  is a semi-algebraic path which
connects the point $y$
to a distinguished point  $a_\theta$
of some roadmap
$\RM(\ZZ({\mathcal P}' \cup \{Q\}, \R^k))$
where ${\mathcal P}'\subset {\mathcal P}$, staying inside
${\RR}(\overline\sigma(y),\ZZ(Q,\R^k)).$
\end{itemize}

Moreover, the complexity of the algorithm is  
$s^{k'+1} d^{O(k^4)}$, where
 $s$ is a bound on the number of elements of
${\mathcal  P}$
and $d$ is a bound on the degrees of $Q$ and
the elements of ${\mathcal  P}$.
\end{theorem}

\subsection{Computing higher Betti numbers}
\label{subsec:nerve}
It is well known that
that 
the Betti numbers of a closed and bounded semi-algebraic set 
can be computed using elementary linear
algebra once we have a triangulation of the set. However, 
triangulations of semi-algebraic sets are expensive
to compute, requiring doubly exponential time. 

One basic idea that underlies some of the recent progress
in designing algorithms for computing the
Betti numbers of semi-algebraic sets is that the cohomology
groups of a semi-algebraic set can often be computed from a sufficiently
well-behaved covering of the set {\em without having to triangulate the set}.

The idea of computing
cohomology from ``good'' covers is an old one in algebraic topology
and the first result in this direction is often called the ``Nerve Lemma''.
In this section we give a brief introduction to the Nerve Lemma and
its generalizations.

We first define formally the notion of a cover of a closed, bounded
semi-algebraic set.
\begin{definition}[Cover]
\label{def:covering}
Let $S \subset \R^k$ be a closed and bounded semi-algebraic set.
A cover, ${\mathcal C}(S)$, of $S$ consists  of an ordered
index set, which by a slight abuse of language we also 
denote by  ${\mathcal C}(S)$, and a map that
associates to each $\alpha \in {\mathcal C}(S)$
a closed and bounded semi-algebraic
subset $S_\alpha \subset S$ such that 
\[
S = \bigcup_{\alpha \in {\mathcal C}(S)} S_\alpha.
\]
\end{definition}

For $\alpha_0,\ldots,\alpha_p, \in {\mathcal C}(S)$,
we associate to the formal product,
$\alpha_0\cdots\alpha_p$, the closed and bounded semi-algebraic set 

\begin{equation}
\label{eqn:defofprod}
S_{\alpha_0\cdots\alpha_p} = S_{\alpha_0} \cap \cdots \cap S_{\alpha_p}.
\end{equation}

Recall that the $0$-th simplicial cohomology group of a closed and bounded 
semi-algebraic set $X$, $\HH^0(X)$, can be identified with the 
$\Q$-vector space
of $\Q$-valued locally constant functions on $X$. It is easy to see that the dimension
of $\HH^0(X)$ is equal to the number of connected components of $X$.

For $\alpha_0,\alpha_1,\ldots,\alpha_p,\beta \in {\mathcal C}(S)$, 
and $\beta \not\in \{\alpha_0,\ldots,\alpha_p\}$,
let
\[
r_{\alpha_0,\ldots,\alpha_p;\beta}: 
\HH^0(S_{\alpha_0\cdots\alpha_p}) \longrightarrow 
\HH^0(S_{\alpha_0\cdots\alpha_{p}\cdot\beta})
\]
be the homomorphism defined as follows.
Given a locally constant function,
$\phi \in \HH^0(S_{\alpha_0\cdots\alpha_p})$,
$r_{\alpha_0\cdots\alpha_p;\beta}(\phi)$ 
is the locally constant function
on $S_{\alpha_0\cdots\alpha_{p}\cdot\beta}$ obtained by restricting
$\phi$ to $S_{\alpha_0\cdots\alpha_{p}\cdot\beta}$.

We define the generalized restriction homomorphisms
\[
\delta^p: \bigoplus_{\alpha_0 < \cdots < \alpha_p, \alpha_i \in {\mathcal C}(S)}  
\HH^0(S_{\alpha_0\cdots\alpha_p})\longrightarrow 
\bigoplus_{\alpha_0< \cdots <\alpha_{p+1}, \alpha_i \in {\mathcal C}(S)} 
\HH^0(S_{\alpha_0\cdots\alpha_{p+1}})
\]
by 
\begin{equation}
\label{eqn:generalized_restriction}
\delta^p(\phi)_{\alpha_0\cdots\alpha_{p+1}}= \sum_{0 \leq i \leq p+1} 
(-1)^{i} r_{\alpha_0\cdots\hat{\alpha_i}\cdots\alpha_{p+1}; \alpha_i}
(\phi_{\alpha_0\cdots\hat{\alpha_{i}}\cdots\alpha_{p+1}}),
\end{equation}
where $\phi \in \bigoplus_{\alpha_0< \cdots <\alpha_p \in {\mathcal C}(S)}  
\HH^0(S_{\alpha_0\cdots\alpha_p})$ and
$r_{\alpha_0\cdots\hat{\alpha_i}\cdots\alpha_{p+1}; \alpha_i}$ 
is the restriction homomorphism defined previously.
The sequence of homomorphisms $\delta^p$ gives rise to a complex,
$\LL^{\bullet}({\mathcal C}(S))$, defined by

\begin{equation}
\label{eqn:defofLL}
\LL^p({\mathcal C}(S)) = 
\bigoplus_{\alpha_0< \cdots <\alpha_p, \alpha_i \in {\mathcal C}(S)} 
\HH^0(S_{\alpha_0\cdots\alpha_p}),
\end{equation}

with the differentials 
$\delta^p: \LL^p({\mathcal C}(S)) \rightarrow \LL^{p+1}({\mathcal C}(S))$ 
defined as in Eqn. (\ref{eqn:generalized_restriction}).

\begin{definition}[Nerve complex]
\label{def:nervecomplex}
The complex $\LL^{\bullet}({\mathcal C}(S))$ is 
called the {\em nerve complex} of the cover  ${\mathcal C}(S)$. 
\end{definition}

For $\ell \geq 0$ we will denote by 
$\LL^{\bullet}_\ell({\mathcal C}(S))$ the truncated complex
defined by

\begin{align*}
\LL^p_\ell({\mathcal C}(S)) &= \LL^p({\mathcal C}(S)),\; 0 \leq p \leq \ell, \\
&= 0,\;  p > \ell.
\end{align*}

Notice that once we have a cover of $S$
and we identify the semi-algebraically connected components of the various intersections,
$S_{\alpha_0\cdots\alpha_p}$, we have natural bases for the vector
spaces
$$
\LL^p({\mathcal C}(S)) = 
\bigoplus_{\alpha_0< \cdots <\alpha_p, \alpha_i \in {\mathcal C}(S)} 
\HH^0(S_{\alpha_0\cdots\alpha_p})
$$ 
appearing as terms of the nerve complex. Moreover, the matrices
corresponding to the homomorphisms $\delta^p$ in this basis depend only
on the inclusion relationships between the semi-algebraically connected components of
$S_{\alpha_0\cdots\alpha_{p+1}}$ and those of 
$S_{\alpha_0\cdots\alpha_p}$. 

\begin{definition}[Leray Property]
\label{def:leray}
We say that the cover
${\mathcal C}(S)$
{\em satisfies the Leray property} if
each non-empty intersection
$S_{\alpha_0\cdots\alpha_p}$ is contractible.
\end{definition}

Clearly, in this case
$$
\begin{array}{cccc}
\HH^0(S_{\alpha_0\cdots\alpha_p}) &\cong & \Q, &\mbox{if $S_{\alpha_0\cdots\alpha_p} \neq \emptyset$} \\
                                  &\cong &  0,  & \mbox{if $S_{\alpha_0\cdots\alpha_p} = \emptyset$}.
\end{array}
$$

It is a classical fact (usually referred to as the {\em Nerve Lemma}) that

\begin{theorem}[Nerve Lemma]
\label{the:nerve}
Suppose that the cover ${\mathcal C}(S)$  satisfies the Leray property.
Then for each $i \geq 0$,
\[ \HH^i(\LL^{\bullet}({\mathcal C}(S))) \cong \HH^i(S).\]
\end{theorem}
(See for instance \cite{Rotman} for a proof.)

\begin{remark}
There are several interesting extensions of Theorem
\ref{the:nerve} (Nerve Lemma). For instance, 
if the Leray property is weakened to say that each $t$-ary intersection
is $(k-t+1)$-connected, then one can conclude that the nerve complex is 
$k$-connected. We refer the reader to the article by Bj\"{o}rner 
\cite{Bjorner} for more details.
\end{remark}

Notice that Theorem \ref{the:nerve} gives a method for 
computing the Betti numbers of $S$ using linear algebra
from a cover of $S$ by contractible sets for which all
non-empty intersections are also contractible, once we are able to  test 
emptiness of the various intersections $S_{\alpha_0\cdots\alpha_p}$.

Now suppose that each individual member, $S_{\alpha_0}$, of the cover is 
contractible, but the various intersections  
$S_{\alpha_0\cdots\alpha_p}$
are not necessarily contractible for $p \geq 1$. 
Theorem \ref{the:nerve} does not hold in this case. 
However, the following theorem is proved in \cite{BPRbettione} and underlies the
singly exponential algorithm for computing the first Betti number of
semi-algebraic sets described there.

\begin{theorem}\cite{BPRbettione}
\label{the:bettione}
Suppose that each individual member, $S_{\alpha_0}$, of the cover 
${\mathcal C}(S)$ is  contractible.
Then,
\[ \HH^i(\LL^{\bullet}_2({\mathcal C}(S))) \cong \HH^i(S),
\mbox{ for } i = 0,1.
\]
\end{theorem}

\begin{remark}
\label{rem:bettione}
Notice that from a cover by contractible sets
Theorem \ref{the:bettione} allows us to compute using linear algebra,
$b_0(S)$ and $b_1(S)$, once we have identified the non-empty 
semi-algebraically connected components
of the pair-wise and triple-wise intersections of the sets in the cover
and their inclusion relationships. 
\end{remark}

\subsubsection{Constructing coverings of closed semi-algebraic sets
by closed contractible sets}
\label{sec:acycliccov}
The parametrized paths obtained in Theorem \ref{the:alg:parametrized} are
not necessarily closed or even contractible, but become so after making
appropriate modifications. At the same time it is possible to maintain the
covering property, namely for any given ${\mathcal P}$-closed semi-algebraic
$S$ set, there exists a set of modified parametrized paths, whose union
is $S$. Moreover, these modified sets are closed and contractible.
We omit the details of this (technical) 
construction  referring the reader to \cite{BPRbettione} for more detail. 
Putting together the constructions outlined above we have:

\begin{theorem}\cite{BPRbettione}
\label{the:alg:covering}
There exists an algorithm that given as input
a  ${\mathcal P}$-closed and bounded semi-algebraic set $S$, 
outputs
a set of formulas $\{\phi_1,\ldots,\phi_M\}$ such that
\begin{itemize}
\item 
each  $\RR(\phi_i,\R'^k)$ is 
semi-algebraically contractible, and
\item 
$\displaystyle{
\bigcup_{1 \leq i \leq M} 
\RR(\phi_i,\R'^k) = \E(S,\R'),
}
$
\end{itemize}
where $\R'$ is some real closed extension of $\R$. 
The complexity of the algorithm is bounded by $s^{(k+1)^2}d^{O(k^5)}$,
where $s = \card \;{\mathcal P}$ and $d = \max_{P \in {\mathcal P}} deg(P).$
\end{theorem}

\subsubsection
{Computing the First Betti Number}
\label{sec:closedcase}
It is now an easy consequence of the existence of singly exponential
time covering algorithm (Theorem \ref{the:alg:covering}), 
and Theorem \ref{the:bettione} stated above, 
along with the fact that we can compute descriptions of the semi-algebraically connected
components of semi-algebraic sets in singly exponential time,
that we can compute the first Betti number of 
closed and bounded semi-algebraic sets in singly exponential time
(see Remark \ref{rem:bettione} above), since the 
dimensions of the images and kernels of
the homomorphisms of the complex,
$\LL^{\bullet}_2({\mathcal C}(S))$ in  Theorem \ref{the:bettione},
can then be computed using traditional algorithms from linear algebra.
As mentioned earlier, 
for arbitrary semi-algebraic sets (not necessarily closed and
bounded),
there is a singly exponential time reduction to the closed and
bounded case using the construction of Gabrielov and Vorobjov
\cite{GaV}.

\subsubsection{Algorithm for Computing the First Few Betti Numbers}
Using the same ideas as above but with a more complicated recursive
procedure to construct a suitable complex one has the following:

\begin{theorem}\cite{Bas05-first}
\label{the:bettifew}
For any given $\ell$, there is an algorithm that
takes as input a ${\mathcal P}$-formula describing a semi-algebraic set 
$S \subset \R^k$,
and outputs $b_0(S),\ldots,b_\ell(S).$ 
The complexity of the algorithm is $(sd)^{k^{O(\ell)}}$, 
where $s  = \card \;({\mathcal P})$ and $d = \max_{P\in {\mathcal P}}{\rm deg}(P).$
\end{theorem}
Note that the complexity is singly exponential in $k$ for every fixed $\ell$.

\subsection{Computing generalized Euler-Poincar\'e characteristic}
\label{sec:eq}
As mentioned before in Section \ref{sec:BKR},
efficient algorithms for sign determination of univariate
polynomials described in \cite{BKR,RS} are amongst
the most basic algorithms in algorithmic real algebraic geometry.
Given
$\mathcal{P} \subset \R[X], Q \in \R[X]$ with
$\card \;\mathcal{P} = s,$ and $\deg(P) \leq d$ for $P \in \mathcal{P} \cup \{Q\}$,
these algorithms count
for each realizable sign condition of the family
$\mathcal{P}$, the cardinality
of the set of real zeros of $Q$, lying in  the realization of that
sign condition. 
The complexity of the algorithm in \cite{RS}
is $sd^{O(1)}$.

In the multidimensional
 case, it is no longer meaningful to talk about the
cardinalities of the zero set of $Q$ lying in the realizations of different
sign conditions of $\mathcal{P}$. However, there exists another discrete
valuation on semi-algebraic sets that properly generalizes  the notion
of cardinality.
This valuation is the Euler-Poincar\'e characteristic.

The \bfdef{Euler-Poincar\'e characteristic}, $\chi(S)$, of
a  closed and bounded semi-algebraic
set $S \subset \R^k$ is defined as
$$\chi(S) = \sum_i (-1)^i b_i(S),$$
where $b_i(S)$ is the rank of the $i$-th simplicial
homology group of $S$. Note that with this definition,
$\chi(\emptyset)=0$, and $\chi(S) = \card \; S,$ whenever
$\card \;S < \infty$. Moreover, $\chi$ is additive. 

The Euler-Poincar\'e characteristic defined above for closed and
bounded semi-algebraic set can be extended additively to all semi-algebraic
sets. This \bfdef{generalized Euler-Poincar\'e characteristic} is then a 
homeomorphism (but not a homotopy) invariant, and establishes an isomorphism
between the \bfdef{Grothendieck ring, $K_0(\mathrm{sa})$},
of homeomorphism classes of semi-algebraic sets
and $\mathbb{Z}$. 

The problem of determining the Euler-Poincar\'e characteristic of 
$\mathcal{P}$-closed
semi-algebraic sets was considered in \cite{Basu1} where an algorithm was
presented for computing the Euler-Poincar\'e characteristic of
a given $\mathcal{P}$-closed semi-algebraic set. 
The complexity of the algorithm is $(ksd)^{O(k)}$. Moreover, in the
special case when the coefficients of the polynomials in $\mathcal{P}$ are
integers of  bit lengths bounded by $\tau$,
the algorithm
performs at most $(ksd)^{O(k)}\tau^{O(1)}$ bit operations.

The following result (which should be viewed as a generalization
of the univariate sign determination algorithm) appears in 
\cite{BPR-euler-poincare}. 

\begin{theorem}
\label{the:alg_BPR-euler-poincare}
There exists an algorithm 
which
given
an algebraic set $Z=\ZZ(Q,\R^k) \subset\R^k$  and a
finite set of polynomials  $\mathcal{P}= \{P_1,\ldots,P_s\}
\subset \R[X_1,\ldots,X_k]$, computes
the list $\chi(\mathcal{P},Z)$ 
indexed by elements, $\sigma,$ of ${\rm Sign}(\mathcal{P},Z)$.
If the degrees of the polynomials in $\mathcal{P} \cup \{Q\}$ are bounded by
$d$, and the real dimension of $Z=\ZZ(Q,\R^k)$ is $k'$, then the
complexity of the algorithm  is
$$s^{k'+1}O(d)^k+ s^{k'}((k'\log_2(s)+k\log_2(d))d)^{O(k)}.$$
If the coefficients of the polynomials in
$\mathcal{P} \cup \{Q\}$ are integers  of bit-sizes bounded by
 $\tau$, then the bit-sizes of the integers
appearing in the intermediate computations and the output
are bounded by
$\tau ((k'\log_2(s)+k\log_2(d))d)^{O(k)}$.
\end{theorem}

\subsection{Relation between the complexity of quantifier elimination and 
the complexity of computing Betti numbers} 
\label{sec:Toda}
It is clear from the previous sections that there are two important strands
of research in algorithms in real algebraic geometry, namely
\begin{enumerate}
\item 
Algorithms for deciding sentences in the first-order theory of the reals
(with several blocks of quantifiers);
\item
Computing topological invariants of semi-algebraic sets (such as their
Betti numbers).
\end{enumerate}

While these two classes of problems might seem quite different, the following 
reduction result gives a polynomial time reduction of the problem of
deciding quantified sentences in the first order theory of the reals with
a fixed number of quantifiers to the problem of computing Betti numbers
of semi-algebraic sets. For technical reasons, the reduction is only proved
for a certain sub-class of formulas which is defined more precisely below. 

\begin{definition}(Compact general decision problem with at most 
$\omega$ quantifier alternations (${\bf GDP_\omega^c}$))

\begin{itemize}
\item[Input.] 
A sentence $\Phi$ in the first order theory of $\R$
\[    (Q_1 \X^{1} \in \Sphere^{k_1})  \cdots 
(Q_\omega \X^{\omega} \in \Sphere^{k_\omega})
\phi(\X^1,\ldots,\X^\omega),
\]
where for each $i, 1 \leq i \leq \omega$,
$\Sphere^{k_i}$ is the unit sphere in $\R^{k_i+1}$,
$\X^i = (X^i_0,\ldots,X^i_{k_i})$ is a block of $k_i+1$ variables,
$Q_i \in \{\exists,\forall\}$, with $Q_j \neq Q_{j+1}, 1 \leq j < \omega$,
and
$\phi$ is a quantifier-free formula defining a {\em closed}
semi-algebraic subset $S$ of $\Sphere^{k_1}\times\cdots\times\Sphere^{k_\omega}$. 
\item[Output.] True or False depending on whether $\Phi$ is true or false
in the first order theory of $\R$.
\end{itemize}
\end{definition}

\begin{notation}
For any semi-algebraic set $S \subset \R^k$, we denote by 
$P_S(T)$, denote the \bfdef{Poincar\'e polynomial} of $S$ -- namely,
\[
P_S(T) := \sum_{i \geq 0} b_i(S)\; T^i.
\]
\end{notation}

\begin{definition}(Computing the Poincar\'e polynomial of 
semi-algebraic sets \\
(\textbf{Poincar\'e}))

\begin{itemize}
\item[Input.] A quantifier-free formula defining a 
semi-algebraic set $S\subset \R^k$.
\item[Output.] 
The Poincar\'e polynomial $P_S(T)$.
\end{itemize}
\end{definition}

The following reduction result appears in \cite{BZ09}. It says
that with a mild hypothesis of compactness, the General Decision
Problem with a fixed number of quantifier alternations can be reduced
in polynomial time to the problem of computing Betti numbers of
semi-algebraic sets.  

\begin{theorem}\cite{BZ09}
\label{the:toda}
For every  $\omega > 0$, there is a deterministic 
polynomial time reduction %%in the  Blum-Shub-Smale model of 
of $\bf{GDP_\omega^c}$ to \upshape \textbf{Poincar\'e}.
\end{theorem}

\begin{remark}
Theorem \ref{the:toda} is motivated by a well known theorem due to
Toda \cite{Toda} in discrete complexity theory which relates 
two complexity classes -- namely the polynomial hierarchy and
the complexity class $\# \mathbf{P}$. Theorem \ref{the:toda} can be viewed
as an analogue of Toda's theorem in the Blum-Shub-Smale model
of computations over arbitrary real closed fields \cite{BCSS98}
(see also \cite{BZ09}).
\end{remark}

The main ingredients in  the proof of Theorem \ref{the:toda} is an
efficient semi-algebraic realization of the iterated fibered join
of a semi-algebraic set with itself over a semi-algebraic map, and
Alexander duality that allows one to express the Poincar\'e polynomial
of a semi-algebraic subset of the sphere in terms of its complement in
the sphere.

\subsection{Effective semi-algebraic triangulation and stratification}
As mentioned above in Section \ref{sec:cad} one obtains an algorithm for
computing a semi-algebraic triangulation of semi-algebraic sets
using cylindrical algebraic decomposition (after making a generic
linear change of co-ordinates). The complexity of this is algorithm 
dominated by the cost of the performing the cylindrical algebraic 
decomposition, and is thus doubly exponential.

Algorithms for computing stratifications of semi-algebraic sets, such
that the strata satisfy additional regularity conditions (such as
Whitney conditions (a) and (b)) have been considered by several authors.
Rannou \cite{Rannou98} gave an algorithm for obtaining stratification
with regularity conditions that imply the Whitney conditions. The complexity
of this algorithm is 
doubly exponential
in the depth of the stratification.
{\em Finding a singly exponential algorithm for computing stratifications of
semi-algebraic sets
remains a major open problem} (see Section \ref{sec:open}).

\subsection{Semi-algebraic sets defined by quadratic and partially quadratic
systems}
\label{sec:quadratic}
A restricted class of semi-algebraic sets - namely, 
semi-algebraic sets defined by quadratic inequalities -- has  been considered
by several researchers \cite{Bar93,Bar97, GrPa04}. 
As in the case of general semi-algebraic sets, the Betti numbers 
of such sets can be exponentially large in the number of variables, 
as can be seen in the following example.

\begin{example}
\label{exa:index}
The set~$S\subset\R^\ell$ defined by 
\[
Y_1(Y_1-1)\ge0,\ldots, Y_\ell(Y_\ell-1)\ge0
\]
satisfies $b_0(S)=2^\ell$.
\end{example}
However, it turns out that for a semi-algebraic set $S \subset \R^{\ell}$
defined by $m$ quadratic inequalities,
it is possible to obtain upper bounds on the Betti numbers of $S$ 
which are polynomial in $\ell$ and exponential only in $m$.
The first such result is due to Barvinok \cite{Bar97}, 
who proved the following  theorem 
(see also \cite{Bas05-first-Kettner,Agrachev-Lerario2012}).

\begin{theorem}\cite{Bar97}
\label{the:barvinok}
Let $S \subset \R^{\ell}$ be defined by $Q_1 \geq 0, \ldots,  Q_m \geq 0$, 
$\deg(Q_i) \leq 2, 1 \leq i \leq m$. Then
$b(S) \leq \ell^{O(m)}$.
\end{theorem}

\begin{remark}
Notice that the bound in Theorem \ref{the:barvinok} is polynomial
in the dimension $\ell$ for fixed $m$, and this fact depends crucially on the 
assumption that the degrees of the polynomials $Q_1,\ldots,Q_m$ are 
at most two.
For instance, the semi-algebraic set defined by a {\em single} 
polynomial of degree $4$ can have Betti numbers exponentially large in 
$\ell$, as exhibited  by the 
semi-algebraic subset of $\R^\ell$ defined by
$$
\displaylines{
\sum_{i=0}^{\ell} Y_i^2(Y_i - 1)^2 \leq 0.
}
$$
The above example illustrates  the delicate nature of the bound in Theorem
\ref{the:barvinok}, since a single inequality of degree $4$ is enough to 
destroy the polynomial nature of the bound. 
In contrast to this, 
it is shown in Theorem \ref{the:BP'R-07} below
that a polynomial bound on the Betti numbers of $S$ continues to hold, 
even if we allow a few (meaning any constant number) of the variables 
to occur with degrees larger than two in the polynomials 
 used to describe the set $S$.
\end{remark}

The bound on the sum of all the Betti numbers in Theorem \ref{the:barvinok}
has exponential dependence on the number of inequalities. This dependence is
unavoidable,  since the semi-algebraic set
$S \subset \R^k$ defined by 
\[
X_1(1 - X_1) \leq 0, \ldots, X_k(1 - X_k) \leq 0,
\]
has $b_0(S) = 2^k$.

Hence, it is somewhat surprising that for any fixed constant $\ell$, 
the Betti numbers
$b_{k-1}(S),\ldots,b_{k-\ell}(S)$, of a basic closed semi-algebraic set
$S \subset \R^k$ defined by quadratic inequalities, are  polynomially 
bounded. The following theorem appears in \cite{B00}.

\begin{theorem}
\label{the:quadratic}
Let  $\R$ a real closed field and $S \subset \R^k$ be defined by 
\[
P_1 \leq 0,\ldots, P_s \leq 0, \deg(P_i) \leq 2, 1 \leq i \leq s.
\]
Then, for $\ell \geq 0$,
\[
b_{k-\ell}(S) \leq {s \choose {\ell}} k^{O(\ell)}.
\]
\end{theorem}

\subsubsection{Algorithm for testing emptiness}
The problem of deciding whether a given semi-algebraic set defined by
a finite set of quadratic inequalities is empty or not was considered
first by Barvinok \cite{Bar93} who proved the following theorem.

\begin{theorem}
\cite{Bar93}
There exists an algorithm which decides if a given system of 
inequalities $Q_1 \geq 0,\ldots, Q_\ell \geq 0$, with each
$Q_i \in \R[X_1,\ldots,X_k], \deg(Q_i) \leq 2$, has a solution
in $\R^k$, whose complexity is bounded by $k^{O(\ell)}$.
\end{theorem}

Barvinok's algorithm did not produce explicit sample points meeting
every semi-algebraically connected component of the set of solutions
(in the style of  Theorem \ref{13:the:samplealg} in the general case). 
This was done by Grigoriev and
Pasechnik \cite{GrPa04}. In fact, they consider the following 
more general situation. 

Let $S \subset \R^k$ be the pull-back of a 
$\mathcal{P}$-semi-algebraic subset $T \subset \R^\ell$ via a 
quadratic map $Q = (Q_1,\ldots,Q_\ell):\R^k \rightarrow \R^\ell$, 
where 
$\mathcal{P} \subset \R[Y_1,\ldots,Y_\ell]$, 
$Q_1,\ldots,Q_\ell \in \R[X_1,\ldots,X_k]$ with $\deg(Q_i) \leq 2$
for $i = 1,\ldots,\ell$.

In \cite{GrPa04},
Grigoriev and Pasechnik give an algorithm that computes a set of sample
points guaranteed to meet every semi-algebraically connected component of $S$ whose complexity
is bounded by $(ksd)^{O(\ell)}$ where $s = \card\; \mathcal{P}$, and
$d$ is a bound on the degrees of the polynomials in $\mathcal{P}$.

\begin{remark}
Note that the problem of deciding the feasibility of even one quartic
real polynomial equation is an $\mathbf{NP}$-hard problem, and the same is 
true for systems of quadratic equations. Thus, there is little
hope for obtaining a polynomial-time algorithm for either of these
problems. The above results are somewhat surprising in that they
imply in the quadratic case one obtains polynomial time algorithms for
testing feasibility, provided the number of polynomials is kept
fixed (see also Section \ref{sec:quadratichardness} below).
We refer the reader to \cite{Pap} and \cite{BCSS98} for precise definitions 
of the computational complexity classes that we refer to here and
elsewhere in this survey.
\end{remark}

\subsubsection{Computing the top few Betti numbers of basic 
semi-algebraic sets defined by quadratic inequalities}
Motivated by the polynomial bound on the top few Betti numbers
of sets defined by quadratic inequalities (Theorem \ref{the:quadratic}),
the problem of obtaining a polynomial time algorithm to compute these
numbers was investigated in \cite{Bas05-top} where the following result
is proved.

\begin{theorem}\cite{Bas05-top}
\label{the:alg:topBetti}
There exists an algorithm 
which given a set of $s$ polynomials,
${\mathcal P} = \{P_1,\ldots,P_s\} \subset \R[X_1,\ldots,X_k],$
with ${\rm deg}(P_i) \leq 2, 1 \leq i \leq s,$
computes $b_{k-1}(S), \ldots, b_{k-\ell}(S),$ 
where $S$ is the set defined by $P_1 \leq 0,\ldots,P_s \leq 0$.
The complexity of the algorithm  is
\begin{equation}
\label{eqn:complexity}
\sum_{i=0}^{\ell+2} {s \choose i} k^{2^{O(\min(\ell,s))}}.
\end{equation}
If the coefficients of the polynomials in
${\mathcal P}$ are integers  of bit-sizes bounded by
 $\tau$, then the bit-sizes of the integers
appearing in the intermediate computations and the output
are bounded by $\tau (sk)^{2^{O(\min(\ell,s))}}.$
\end{theorem}

\subsubsection{Significance from the computational 
complexity theory viewpoint}
\label{sec:quadratichardness}
Semi-algeb\-raic sets defined by a system of quadratic inequalities 
have a special significance in the theory of computational complexity. 
Even though such sets might seem to be the next simplest class of 
semi-algebraic sets after sets defined by linear inequalities, 
from the point of view of 
computational complexity they represent a quantum leap.
Whereas there exist (weakly) polynomial time algorithms for
solving linear programming, solving quadratic
feasibility problem is provably hard.
For instance, it follows from an easy reduction from the problem of 
testing feasibility of a real quartic equation in many variables, that 
the problem of testing whether a system of quadratic inequalities 
is feasible is  $\mathbf{NP}_{\rm R}$-complete in the Blum-Shub-Smale
model of computation (see \cite{BCSS98}). 
Assuming the input polynomials to have integer
coefficients, the same problem is $\mathbf{NP}$-hard 
in the classical Turing machine
model, since it is also not difficult to see that the Boolean satisfiability 
problem can be posed as the problem of deciding whether a certain 
semi-algebraic set  defined by  quadratic inequalities is empty or not.
Counting the number of semi-algebraically connected 
components of such sets is even harder. In fact, it is 
$\mathbf{PSPACE}$-hard \cite{Reif79}
($\mathbf{PSPACE}$ is a complexity class which
contains the entire polynomial hierarchy),
and the proof of this
results extend easily to the  quadratic case.  
Moreover, it is proved in \cite{Bas05-top}
for $\ell = O(\log k)$,
computing the $\ell$-th Betti number of 
a basic semi-algebraic set defined by quadratic inequalities in $\R^k$
is $\mathbf{PSPACE}$-hard.
In view of these hardness results,  it is unlikely that there exist polynomial
time algorithms for computing the Betti numbers (or even the first
few Betti numbers) of such a set.

From this point of view, 
Theorem \ref{the:alg:topBetti} is quite surprising, since it gives a 
polynomial time algorithm for computing certain Betti numbers of 
a class of semi-algebraic sets for which computing the zero-th Betti
number is already $\mathbf{PSPACE}$-hard.

\subsubsection{Semi-algebraic sets defined by partially quadratic systems}
We have discussed topological as well as algorithmic results concerning
general semi-algebraic sets, as well as those defined by quadratic
constraints
In \cite{BP'R07jems}, the authors try to interpolate between 
results known for general  semi-algebraic sets (defined by polynomials
of arbitrary degrees) and those known for semi-algebraic sets defined
by polynomials of degree at most $2$.
In order to do so they consider  semi-algebraic sets defined by 
polynomial inequalities,
in which the dependence of the 
polynomials on a {\em subset of the variables}  is at most quadratic.
As a result we obtain common generalizations of the bounds stated in 
Theorems \ref{the:B99} and \ref{the:barvinok}.
Given any polynomial $P \in \R[X_1,\ldots,X_k,Y_1,\ldots,Y_\ell]$, we will
denote by $\deg_X(P)$ (resp. $\deg_Y(P)$) the total degree of $P$ with respect
to the variables $X_1,\ldots,X_k$ (resp. $Y_1,\ldots,Y_\ell$).

Denote by 
\begin{itemize}
\item
${\mathcal Q}\subset  \R[Y_1,\ldots,Y_\ell,X_1,\ldots,X_k]$,
a family of polynomials
with 
\[
\deg_{Y}(Q) \leq 2, 
\deg_{X}(Q) \leq d,  Q\in {\mathcal Q}, \card \;{\mathcal Q}=m,
\]
\item
${\mathcal P} \subset \R[X_1,\ldots,X_k]$,
a family of polynomials 
with
\[
\deg_{X}(P) \leq d, P \in {\mathcal P}, \card \;{\mathcal P}=s.
\]
\end{itemize}

The following theorem that interpolates between 
Theorems \ref{the:GV} and \ref{the:barvinok} 
above is proved in \cite{BP'R07jems}.

\begin{theorem}
\label{the:BP'R-07}
Let $S \subset \R^{\ell+k}$ 
be a $({\mathcal P} \cup {\mathcal Q})$-closed semi-algebraic set. Then
$$
\displaylines{
b(S) \leq 
\ell^2 (O(s+\ell+m)\ell d)^{k+2m}. 
}
$$
In particular, for $m \leq \ell$, we have
$
\displaystyle{
b(S) \leq \ell^2 (O(s+\ell)\ell d)^{k+2m}. 
}
$
\end{theorem}

Notice that Theorem \ref{the:BP'R-07} 
can be seen as a common generalization of 
Theorems \ref{the:B99} and \ref{the:barvinok},
in the sense that 
we recover similar bounds (that is bounds having the same shape) 
as in Theorem \ref{the:B99} (respectively Theorem \ref{the:barvinok})
by setting $\ell$ and $m$ 
(respectively $s$, $d$ and $k$) to $O(1)$.

Note also that as a special case of Theorem \ref{the:BP'R-07}
we obtain a bound on the sum of the Betti numbers of
a semi-algebraic set defined over a quadratic map.
As mentioned before,
such sets have been considered from an algorithmic point of view
in \cite{GrPa04}, where 
an efficient algorithm is described for computing sample 
points in every semi-algebraically connected component, as well as testing emptiness, 
of such sets.

More precisely we have: 
\begin{corollary}
\label{cor:main}
Let $Q = (Q_1,\ldots,Q_k) : \R^{\ell} \rightarrow \R^k$ be a map where
each $Q_i \in \R[Y_1,\ldots,Y_\ell]$ and $\deg(Q_i) \leq 2$.
Let $V \subset \R^k$ be a ${\mathcal P}$-closed semi-algebraic set for
some family ${\mathcal P} \subset \R[X_1,\ldots,X_k]$, with
$\card \;{\mathcal P} = s$ and $\deg(P) \leq d, P \in {\mathcal P}$.
Let $S = Q^{-1}(V)$. Then
$$
\displaylines{
b(S) \leq \ell^2 (O(s+\ell+k)\ell d)^{3k}. 
}
$$
\end{corollary}

The techniques developed in this paper for obtaining tight bounds on the
Betti numbers of semi-algebraic sets defined by partly quadratic systems of
polynomials also pave the way towards designing more 
efficient algorithms for computing the Euler-Poincar\'e characteristic as
well as the Betti numbers of such sets. 

The following theorem appears in \cite{BP'R07jems}.
\begin{theorem}
\label{the:algo-EP}
There exists an algorithm 
that takes as input the description of a
$({\mathcal P} \cup {\mathcal Q})$-closed semi-algebraic set
$S$ (following the same notation as in Theorem \ref{the:BP'R-07})
and outputs its
the Euler-Poincar\'e characteristic
$\chi(S)$. 
The complexity of this algorithm  is bounded by 
$(\ell s m d)^{O(m(m+k))}$. 
In the case when $S$ is a basic closed semi-algebraic set
the complexity of the algorithm is 
$(\ell s m d)^{O(m+k)}$.
\end{theorem}

The algorithm for computing all the Betti numbers has complexity
$(\ell s m d)^{2^{O(m+k)}}$ and 
its description can be found in \cite{BP'R07joa}. 
While the complexity of both the algorithms discussed above 
is {\em polynomial}  for fixed $m$ and $k$, the complexity of the algorithm
for computing the Euler-Poincar\'e characteristic is significantly better
than that of the algorithm for computing all the Betti numbers. 

Note that the first versions of both these algorithms for computing
the Euler-Poincar\'e characteristic as well as the 
Betti numbers
of semi-algebraic sets defined by purely quadratic constraints
having complexity which is polynomial for fixed number of constraints,
appeared first in \cite{Bas05-euler} and \cite{Bas05-top} respectively.
The extensions of these algorithms to semi-algebraic sets defined by
partially quadratic systems were made in \cite{BP'R07jems} and 
\cite{BP'R07joa} respectively.

These latter  results indicate that 
the problem of computing the
Betti numbers of semi-algebraic sets defined by a constant number of
polynomial inequalities is solvable in polynomial time, even if we allow
a small (constant sized) subset of the variables to occur with degrees
larger than two in the polynomials defining the given set.

\section{Sums of squares and semi-definite programming}
\label{sec:sos}
All the algorithms surveyed above have the feature that they are exact,
and most of them work over arbitrary real closed fields (even non-archimedean
ones). For example, the
ring generated by the coefficients, $\D$,  
could be the ordered ring, $\Z[\eps]$ with $\eps$
positive and infinitesimal, contained in the real closed field
$\R = \re_{\mathrm{alg}}\langle\eps\rangle$ and all algorithms reported
above would still work without any modification.

There are some other approaches to designing algorithms for
solving systems of real polynomial equations or testing emptiness of
semi-algebraic sets that deserve mention. These approaches strictly
assume that the underlying real closed field is the field $\re$ of real
numbers, and the computations are done with some finite precision. In other
words, the algorithms are numerical rather than exact, and as such there is
some possibility of error in the outputs. These algorithms
are often used in practical applications, where exact or symbolic algorithms
are deemed to be too expensive and small errors are considered not very
significant.

We mention one such approach below. 

\subsection{Deciding non-negativity of polynomials using sums-of-squares}
The problem is to decide whether a given polynomial $P \in \re[X_1,\ldots,X_k]$
is non-negative in $\re^k$. More generally, the problem is to decide
whether a given polynomial $P \in \re[X_1,\ldots,X_k]$ 
is non-negative over a given basic, 
semi-algebraic subset $K \subset \re^k$. 

There are also optimization versions of these problems namely. 

Given $P \in \re[X_1,\ldots,X_k]$ compute 
\[
p^{\min} := \inf_{x \in \re^k} P(x).
\]

More generally,
Given $P \in \re[X_1,\ldots,X_k]$ and $K \subset \re^k$ a basic semi-algebraic
set, compute 
\[
p^{\min} := \inf_{x \in K} P(x).
\]

For purposes of exposition we 
concentrate on the first versions of these problems.

Let the degree of $P$ be $2d$ and
let $\Pos_{k,d}$ (resp. $\Sigma_{k,d}$) denote the cone of non-negative
polynomials (resp. cone of sum of squares) in $\re[X_1,\ldots,X_k]$ of
degree at most $2d$. Clearly $\Sigma_{k,d} \subset \Pos_{k,d}$ and 
as known since Hilbert,
the
inclusion is strict unless the pair $(k,d)$ is of the form 
$(1,d), (k,1)$ or  $(k,d) = (2,2)$ \cite[Chapter 6]{BCR}.
Note that the cones $\Pos_{k,d}$ are in general not understood very well
(for instance, their face structure, extreme rays etc.) and testing membership
in them is clearly an $\mathbf{NP}$-hard problem. 
On the other hand, the cones
$\Sigma_{k,d}$ are relatively well understood and membership in
$\Sigma_{k,d}$ can be tested via semi-definite programming
as a result of the following theorem.

For any symmetric, square matrix $X \in \re^{k \times k}$, we let
$X \succeq 0$ denote that $X$ is positive, semi-definite.
For each $k,d \geq 0$, we denote by $\mathcal{M}_{k,d}$ the set of
exponent vectors $\alpha = (\alpha_1,\ldots,\alpha_k) 
\in \mathbb{N}^k$ with 
$|\alpha| = \sum_{i=1}^{k} \alpha_i \leq d$.
 
\begin {theorem}\cite{Choi-Lam-Reznick92,Powers-Wormann98} 
\label{thm:sdp}
The following are equivalent.
\begin{enumerate}
\item
$P = \sum_{\alpha \in \mathcal{M}_{k,d}} p_\alpha X^\alpha
\in \Sigma_{k,d}$.
\item
The following system in matrix variables 
$X = (X_{\alpha,\beta})_{\alpha,\beta \in \mathcal{M}_{k,d}}$  is feasible:
$$
\displaylines{
X \succeq  0 \cr
\sum_{\beta,\gamma \in \mathcal{M}_{k,d}, \beta+\gamma = \alpha} 
X_{\beta,\gamma} = p_\alpha, \alpha \in \mathcal{M}_{k,2d}.
}
$$
\end{enumerate}
\end{theorem}

The feasibility problem in the above theorem is an instance of 
the feasibility problem in the theory of \bfdef{semi-definite programming}. 
Semi-definite programming (or semi-definite optimization)
is a generalization of linear programming, 
where the problem is to optimize a linear functional over some affine
section of the cone of real symmetric positive semi-definite
matrices in the space of $k \times k$ real symmetric matrices.
Because of its wide ranging applicability, semi-definite programming
has been the focus of  intense effort on the part of researchers 
in optimization for developing efficient algorithms for solving 
semi-definite programming problems. As a result
very efficient algorithms based on ``interior point methods'' 
(see \cite{NN93})
have been developed 
for solving
semi-definite optimization problems such as the one in 
Theorem \ref{thm:sdp}.
These algorithms are very efficient in practice, but there 
seems to be no definitive
mathematical result which states that the running time is polynomial (in the
bit-size of the input) (unlike in the case of linear programming).

Note that the polynomial optimization problems can also be ``approximated''
using the sum of squares cone just like above. For example, in order to 
compute
\[
p^{\min} := \inf_{x \in \re^k} P(x) = 
\sup \{\rho \in \re\;  \mid \;P-\rho \in \Pos_{k,d} \},
\]
one computes
\[
p^{\mathrm{sos}} := \sup \{\rho \in \re\;  \mid \;P-\rho \in \Sigma_{k,d} \}.
\]

Since, 
this latter problem is an example of semi-definite {\em optimization}
problem and can be solved in practice using efficient interior points 
methods. Also note that since the latter problem involves optimization
over a smaller cone we  have that
\[
p^{\mathrm{sos}} \leq p^{\min}.
\]

The idea of ``relaxing''  polynomial optimization problems to 
semi-definite programming has been utilized by Lasserre \cite{Lasserre01,
Lasserre06-1,Lasserre06-2}, Parrilo \cite{Parrilo03} and others to
obtain algorithms for performing polynomial optimization
which perform well in practice (but see Remark \ref{rem:blekherman} below).
 
\begin{remark}
\label{rem:blekherman}
While the idea of approximating the cone of non-negative polynomials by the
smaller cone of sums of square seems to work well in practice for solving
or approximating well solutions of polynomial optimization problems,
one should be aware of certain negative results. Blekherman 
\cite{Blekherman06} proved that the ratio of the volumes of certain
fixed sections of the  cones
$\Sigma_{k,d}$ and $\Pos_{k,d}$ goes to $0$ with $k$ exponentially fast.   
This seems to indicate that the approximation of $\Pos_{k,d}$ by
$\Sigma_{k,d}$ is very inaccurate as $k$ grows (with $d$ fixed).
\end{remark}

We refer the reader to the excellent survey article by Laurent 
\cite{Laurent09} for
more detailed information about the sums-of-squares methods in algorithmic
real algebraic geometry.

\subsection{Complexity of semi-definite programming}
Since semi-definite optimization problems play an important role in the
sums-of-square approximation algorithms described above, it is important
to be aware of the current complexity status of this problem.
As noted above, while interior points algorithms for solving semi-definite
programming problems are extremely efficient in practice, there is no
definite result known placing the semi-definite
programming problem in the class $\mathbf{P}$.
Khachiyan and Prokolab \cite{KP97} proved that there exists a polynomial
time algorithm for semi-definite
programming in case the dimension is fixed. 
Using results proved by Ramana \cite{Ramana97} on exact semi-definite
duality theory, it can be deduced (see \cite{Tarasov-Vyalyi08}) 
that semi-definite feasibility cannot be $\mathbf{NP}$-complete
unless $\mathbf{NP} = \mathbf{co\mbox{-}NP}$ (a hypothesis not believed to be
true). In the Blum-Shub-Smale model of computation over real machines
\cite{BCSS98},
the semi-definite feasibility  problem is clearly in the class 
$\mathbf{NP}_{\re}$, 
and it is unknown
if it is any 
easier than ordinary real polynomial feasibility problem in this 
model.

\section{Open problems}
\label{sec:open}
We list here some interesting open problems some of which could possibly
be tackled in the near future. \\

\paragraph{\em Computing Betti numbers in singly exponential time?}

Suppose $S \subset \R^k$ is a semi-algebraic set defined in terms
of $s$ polynomials, of degrees bounded by $d$. One of the most fundamental
open questions in algorithmic semi-algebraic geometry, is whether
there exists a singly exponential (in $k$) time algorithm for computing
the Betti numbers of $S$. 
The best we can do so far is summarized in Theorem \ref{the:bettifew}
which gives the existence of singly exponential time algorithms for
computing the first $\ell$ Betti numbers of $S$ for any constant
$\ell$. 
A big challenge is 
to design an algorithm
for computing all the Betti numbers of $S$. \\

\paragraph{\em Computing semi-algebraic triangulations in singly
exponential time?}
A related question is whether 
there exists an algorithm for 
computing  semi-algebraic triangulations
with singly exponential complexity. Clearly, such an algorithm
would also make possible the computation of Betti numbers
in singly exponential time. \\

\paragraph{\em More Efficient Algorithms for 
Computing the Number of Connected Components in the Quadratic
Case?}

As described in Section \ref{sec:quadratic}
for semi-algebraic sets in $\R^k$ defined by $\ell$ quadratic inequalities,
there are algorithms for deciding emptiness, as well as computing
sample points in every semi-algebraically connected component whose complexity is
bounded by $k^{O(\ell)}$.
We also have an algorithm
for computing the Euler-Poincar\'e characteristic of such sets whose
complexity is  $k^{O(\ell)}$. 
However, the best known algorithm for computing
the number of semi-algebraically connected components of such sets 
has complexity $k^{2^{O(\ell)}}$ (as a special case of the 
algorithm for computing all the Betti numbers given in Theorem
\ref{the:quadratic}). This raises the question whether there
exists a more efficient algorithm with complexity $k^{O(\ell)}$ or even
$k^{O(\ell^2)}$ for counting the number of semi-algebraically connected components of such sets.
Roadmap type constructions
used for counting semi-algebraically connected components in the case of 
general semi-algebraic sets cannot be directly employed in this context,
because such algorithms will have complexity exponential in $k$. 
Recent work by Coste and Moussa \cite{Coste-Moussa10} on the geodesic
diameter of semi-algebraic sets defined by few quadratic inequalities
might contain some relevant hints towards this goal.
\\

\paragraph{\em More Efficient Algorithms for 
Computing the Number of Connected Components for General Semi-algebraic Sets?}
A very interesting open question is whether the exponent $O(k^2)$ 
in the complexity of roadmap algorithms 
(cf. Theorem \ref{16:the:saconnecting}) can be improved to $O(k)$, so that the
complexity of testing connectivity becomes asymptotically the same as that
of testing emptiness of a semi-algebraic set (cf. Theorem 
\ref{13:the:samplealg}). Recent improvements in the complexity
of roadmap algorithms described in Section \ref{sec:baby-giant} above, 
certainly gives some hope
in this regard.

Such an improvement
would go a long way in making this algorithm practically useful.
It would also be of interest for studying 
metric properties of semi-algebraic sets because of the following.
Applying Crofton's formula from integral geometry  one immediately
obtains as a corollary of Theorem \ref{16:the:saconnecting}
(using the same notation as in the theorem) 
an upper bound of $s^{k'+1}d^{O(k^2)}$
on the length of a semi-algebraic connecting
path connecting two points in any semi-algebraically connected component of $S$
(assuming that $S$ is contained in the unit ball centered at the origin).
An improvement in the complexity of algorithms for constructing
connecting paths (such as the roadmap algorithm) 
would also improve the bound on the length of connecting paths.
Recent results due to D'Acunto and  Kurdyka \cite{DK} show that 
it is possible to construct  semi-algebraic paths
of length $d^{O(k)}$ between two points of $S$
(assuming that $S$ is a semi-algebraically connected
component of a real algebraic set contained in the  unit ball
defined by polynomials of degree $d$). 
However, the semi-algebraic complexity
of such paths cannot be bounded in terms of the parameters $d$ and $k$.
The improvement in the complexity suggested above, apart from
its algorithmic significance,  would also be an
effective version of the results in \cite{DK}. \\

\paragraph{\em Remove the compactness assumption in Theorem
\ref{the:toda}.} More generally, investigate the role of
compactness in the Blum-Shub-Smale model of computations over
real closed fields (see \cite{BZ09} for more details). \\

\paragraph{\em Studying complexity of constructible sheaves and functors} Very recently,
in \cite{Basu2014} a study of \emph{constructible sheaves and functors} have been undertaken from
a complexity point of view. The closure of the category of constructible sheaves under certain operations
can be seen as a far reaching generalization of the Tarski-Seidenberg principle, with many applications. 
Many of the techniques and algorithms developed for semi-algebraic sets and described in this
survey should prove useful in this more general context. 

\bibliographystyle{abbrv}
\bibliography{master}

\def\cprime{$'$}
\begin{thebibliography}{10}

\bibitem{Agrachev-Lerario2012}
A.~Agrachev and A.~Lerario.
\newblock Systems of quadratic inequalities.
\newblock {\em Proc. London Math. Soc.(3)}, 2012.

\bibitem{BGHM97}
B.~Bank, M.~Giusti, J.~Heintz, and G.~M. Mbakop.
\newblock Polar varieties, real equation solving, and data structures: the
  hypersurface case.
\newblock {\em J. Complexity}, 13(1):5--27, 1997.

\bibitem{BGHM01}
B.~Bank, M.~Giusti, J.~Heintz, and G.~M. Mbakop.
\newblock Polar varieties and efficient real elimination.
\newblock {\em Math. Z.}, 238(1):115--144, 2001.

\bibitem{BGHMS10}
B.~Bank, M.~Giusti, J.~Heintz, M.~Safey El~Din, and {\'E}.~Schost.
\newblock On the geometry of polar varieties.
\newblock {\em Appl. Algebra Engrg. Comm. Comput.}, 21(1):33--83, 2010.

\bibitem{Barone-Basu11a}
S.~Barone and S.~Basu.
\newblock Refined bounds on the number of connected components of sign
  conditions on a variety.
\newblock {\em Discrete Comput. Geom.}, 2011.

\bibitem{Bar93}
A.~I. Barvinok.
\newblock Feasibility testing for systems of real quadratic equations.
\newblock {\em Discrete Comput. Geom.}, 10(1):1--13, 1993.

\bibitem{Bar97}
A.~I. Barvinok.
\newblock On the {B}etti numbers of semialgebraic sets defined by few quadratic
  inequalities.
\newblock {\em Math. Z.}, 225(2):231--244, 1997.

\bibitem{B99b}
S.~Basu.
\newblock New results on quantifier elimination over real closed fields and
  applications to constraint databases.
\newblock {\em J. ACM}, 46(4):537--555, 1999.

\bibitem{Basu1}
S.~Basu.
\newblock On bounding the {B}etti numbers and computing the {E}uler
  characteristic of semi-algebraic sets.
\newblock {\em Discrete Comput. Geom.}, 22(1):1--18, 1999.

\bibitem{B00}
S.~Basu.
\newblock Different bounds on the different {B}etti numbers of semi-algebraic
  sets.
\newblock {\em Discrete Comput. Geom.}, 30(1):65--85, 2003.
\newblock ACM Symposium on Computational Geometry (Medford, MA, 2001).

\bibitem{Bas05-first}
S.~Basu.
\newblock Computing the first few {B}etti numbers of semi-algebraic sets in
  single exponential time.
\newblock {\em J. Symbolic Comput.}, 41(10):1125--1154, 2006.

\bibitem{Bas05-euler}
S.~Basu.
\newblock Efficient algorithm for computing the {E}uler-{P}oincar\'e
  characteristic of a semi-algebraic set defined by few quadratic inequalities.
\newblock {\em Comput. Complexity}, 15(3):236--251, 2006.

\bibitem{Bas05-top}
S.~Basu.
\newblock Computing the top few {B}etti numbers of semi-algebraic sets defined
  by quadratic inequalities in polynomial time.
\newblock {\em Found. Comput. Math.}, 8(1):45--80, 2008.

\bibitem{Basu2014}
S.~{Basu}.
\newblock {A complexity theory of constructible functions and sheaves}.
\newblock {\em ArXiv e-prints}, Sept. 2013.

\bibitem{BGV2013}
S.~Basu, A.~Gabrielov, and N.~Vorobjov.
\newblock Monotone functions and maps.
\newblock {\em Rev. R. Acad. Cienc. Exactas F\'\i s. Nat. Ser. A Math. RACSAM},
  107(1):5--33, 2013.

\bibitem{Bas05-first-Kettner}
S.~Basu and M.~Kettner.
\newblock A sharper estimate on the {B}etti numbers of sets defined by
  quadratic inequalities.
\newblock {\em Discrete Comput. Geom.}, 39(4):734--746, 2008.

\bibitem{BP'R07joa}
S.~Basu, D.~V. Pasechnik, and M.-F. Roy.
\newblock Computing the {B}etti numbers of semi-algebraic sets defined by
  partly quadratic sytems of polynomials.
\newblock {\em J. Algebra}, 321(8):2206--2229, 2009.

\bibitem{BP'R07jems}
S.~Basu, D.~V. Pasechnik, and M.-F. Roy.
\newblock Bounding the {B}etti numbers and computing the {E}uler-{P}oincar\'e
  characteristic of semi-algebraic sets defined by partly quadratic systems of
  polynomials.
\newblock {\em J. Eur. Math. Soc. (JEMS)}, 12(2):529--553, 2010.

\bibitem{BPR8}
S.~Basu, R.~Pollack, and M.-F. M.-F.~Roy.
\newblock On the {B}etti numbers of sign conditions.
\newblock {\em Proc. Amer. Math. Soc.}, 133(4):965--974 (electronic), 2005.

\bibitem{BPR95}
S.~Basu, R.~Pollack, and M.-F. Roy.
\newblock On the combinatorial and algebraic complexity of quantifier
  elimination.
\newblock {\em J. ACM}, 43(6):1002--1045, 1996.

\bibitem{BPR95b}
S.~Basu, R.~Pollack, and M.-F. Roy.
\newblock On computing a set of points meeting every cell defined by a family
  of polynomials on a variety.
\newblock {\em J. Complexity}, 13(1):28--37, 1997.

\bibitem{BPR99}
S.~Basu, R.~Pollack, and M.-F. Roy.
\newblock Computing roadmaps of semi-algebraic sets on a variety.
\newblock {\em J. Amer. Math. Soc.}, 13(1):55--82, 2000.

\bibitem{BPR10}
S.~Basu, R.~Pollack, and M.-F. Roy.
\newblock Betti number bounds, applications and algorithms.
\newblock In {\em Current Trends in Combinatorial and Computational Geometry:
  Papers from the Special Program at MSRI}, volume~52 of {\em MSRI
  Publications}, pages 87--97. Cambridge University Press, 2005.

\bibitem{BPR-euler-poincare}
S.~Basu, R.~Pollack, and M.-F. Roy.
\newblock Computing the {E}uler-{P}oincar\'e characteristics of sign
  conditions.
\newblock {\em Comput. Complexity}, 14(1):53--71, 2005.

\bibitem{BPRbook2}
S.~Basu, R.~Pollack, and M.-F. Roy.
\newblock {\em Algorithms in real algebraic geometry}, volume~10 of {\em
  Algorithms and Computation in Mathematics}.
\newblock Springer-Verlag, Berlin, 2006 (second edition).

\bibitem{BPRbettione}
S.~Basu, R.~Pollack, and M.-F. Roy.
\newblock Computing the first {B}etti number of a semi-algebraic set.
\newblock {\em Found. Comput. Math.}, 8(1):97--136, 2008.

\bibitem{BR10}
S.~Basu and M.-F. Roy.
\newblock Bounding the radii of balls meeting every connected component of
  semi-algebraic sets.
\newblock {\em J. Symbolic Comput.}, 45(12):1270--1279, 2010.

\bibitem{BR14}
S.~Basu and M.-F. Roy.
\newblock Divide and conquer roadmap for algebraic sets.
\newblock {\em Discrete Comput. Geom.}, 52(2):278--343, 2014.

\bibitem{BRMS10}
S.~{Basu}, M.-F. {Roy}, M.~{Safey El Din}, and {\'E}.~{Schost}.
\newblock {A baby step-giant step roadmap algorithm for general algebraic
  sets}.
\newblock {\em ArXiv e-prints}, Jan. 2012.

\bibitem{BV06}
S.~Basu and N.~Vorobjov.
\newblock On the number of homotopy types of fibres of a definable map.
\newblock {\em J. Lond. Math. Soc. (2)}, 76(3):757--776, 2007.

\bibitem{BZ09}
S.~Basu and T.~Zell.
\newblock Polynomial hierarchy, {B}etti numbers, and a real analogue of
  {T}oda's theorem.
\newblock {\em Found. Comput. Math.}, 10(4):429--454, 2010.

\bibitem{BKR}
M.~Ben-Or, D.~Kozen, and J.~Reif.
\newblock The complexity of elementary algebra and geometry.
\newblock {\em J. of Computer and Systems Sciences}, 18:251--264, 1986.

\bibitem{Libkin98}
M.~Benedikt, G.~Dong, L.~Libkin, and L.~Wong.
\newblock Relational expressive power of constraint query languages.
\newblock {\em J. ACM}, 45(1):1--34, 1998.

\bibitem{Bjorner}
A.~Bj\"{o}rner.
\newblock Topological methods.
\newblock In R.~Graham, M.~Grotschel, and L.~Lovasz, editors, {\em Handbook of
  Combinatorics}, volume~II, pages 1819--1872. North-Holland/Elsevier, 1995.

\bibitem{Blekherman06}
G.~Blekherman.
\newblock There are significantly more nonnegative polynomials than sums of
  squares.
\newblock {\em Israel J. Math.}, 153:355--380, 2006.

\bibitem{BCSS98}
L.~Blum, F.~Cucker, M.~Shub, and S.~Smale.
\newblock {\em Complexity and real computation}.
\newblock Springer-Verlag, New York, 1998.
\newblock With a foreword by Richard M. Karp.

\bibitem{BCR}
J.~Bochnak, M.~Coste, and M.-F. Roy.
\newblock {\em G\'eom\'etrie alg\'ebrique r\'eelle (Second edition in english:
  Real Algebraic Geometry)}, volume 12 (36) of {\em Ergebnisse der Mathematik
  und ihrer Grenzgebiete [Results in Mathematics and Related Areas ]}.
\newblock Springer-Verlag, Berlin, 1987 (1998).

\bibitem{BoCaR}
F.~Boudaoud, F.~Caruso, and M.-F. Roy.
\newblock Certificates of positivity in the {B}ernstein basis.
\newblock {\em Discrete Comput. Geom.}, 39(4):639--655, 2008.

\bibitem{Canny93a}
J.~Canny.
\newblock Computing road maps in general semi-algebraic sets.
\newblock {\em The Computer Journal}, 36:504--514, 1993.

\bibitem{Choi-Lam-Reznick92}
M.~D. Choi, T.~Y. Lam, and B.~Reznick.
\newblock Sums of squares of real polynomials.
\newblock In {\em {$K$}-theory and algebraic geometry: connections with
  quadratic forms and division algebras ({S}anta {B}arbara, {CA}, 1992)},
  volume~58 of {\em Proc. Sympos. Pure Math.}, pages 103--126. Amer. Math.
  Soc., Providence, RI, 1995.

\bibitem{Collins75}
G.~Collins.
\newblock {Quantifier elimination for real closed fields by cylindrical
  algebraic decomposition}.
\newblock {\em LNCS}, 33:134--183, 1975.

\bibitem{Col}
G.~E. Collins.
\newblock Quantifier elimination for real closed fields by cylindric algebraic
  decomposition.
\newblock In {\em Second GI Conference on Automata Theory and Formal
  Languages}, volume~33 of {\em Lecture Notes in Computer Science}, pages
  134--183, Berlin, 1975. Springer- Verlag.

\bibitem{Coste-Moussa10}
M.~Coste and S.~Moussa.
\newblock Geodesic diameter of sets defined by few quadratic equations and
  inequalities.
\newblock {\em arXiv:1004.5047v1 [math.LO]}.

\bibitem{CR}
M.~Coste and M.-F. Roy.
\newblock Thom's lemma, the coding of real algebraic numbers and the topology
  of semi-algebraic sets.
\newblock {\em Journal of Symbolic Computation}, 5(1/2):121--129, 1988.

\bibitem{DK}
D.~D'Acunto and K.~Kurdyka.
\newblock Bounds for gradient trajectories and geodesic diameter of real
  algebraic sets.
\newblock {\em Bull. London Math. Soc.}, 38(6):951--965, 2006.

\bibitem{Davenport-Heintz88}
J.~H. Davenport and J.~Heintz.
\newblock Real quantifier elimination is doubly exponential.
\newblock {\em J. Symbolic Comput.}, 5(1-2):29--35, 1988.

\bibitem{GaV}
A.~Gabrielov and N.~Vorobjov.
\newblock Betti numbers of semialgebraic sets defined by quantifier-free
  formulae.
\newblock {\em Discrete Comput. Geom.}, 33(3):395--401, 2005.

\bibitem{GV07}
A.~Gabrielov and N.~Vorobjov.
\newblock Approximation of definable sets by compact families, and upper bounds
  on homotopy and homology.
\newblock {\em J. Lond. Math. Soc. (2)}, 80(1):35--54, 2009.

\bibitem{GR92}
L.~Gournay and J.~J. Risler.
\newblock Construction of roadmaps of semi-algebraic sets.
\newblock {\em Appl. Algebra Eng. Commun. Comput.}, 4(4):239--252, 1993.

\bibitem{GrPa04}
D.~Grigoriev and D.~V. Pasechnik.
\newblock Polynomial-time computing over quadratic maps. {I}. {S}ampling in
  real algebraic sets.
\newblock {\em Comput. Complexity}, 14(1):20--52, 2005.

\bibitem{GV92}
D.~Grigoriev and N.~Vorobjov.
\newblock Counting connected components of a semi-algebraic set in
  subexponential time.
\newblock {\em Comput. Complexity}, 2(2):133--186, 1992.

\bibitem{GV}
D.~Y. Grigoriev and N.~N. Vorobjov, Jr.
\newblock Solving systems of polynomial inequalities in subexponential time.
\newblock {\em J. Symbolic Comput.}, 5(1-2):37--64, 1988.

\bibitem{Hardt}
R.~Hardt.
\newblock Semi-algebraic local-triviality in semi-algebraic mappings.
\newblock {\em Amer. J. Math.}, 102(2):291--302, 1980.

\bibitem{HRS94}
J.~Heintz, M.-F. Roy, and P.~Solern\`{o}.
\newblock Description of the connected components of a semialgebraic set in
  single exponential time.
\newblock {\em Discrete and Computational Geometry}, 11:121--140, 1994.

\bibitem{JP2011}
G.~Jeronimo and D.~Perrucci.
\newblock On the minimum of a positive polynomial over the standard simplex.
\newblock {\em J. Symb. Comput.}, 45(4):434--442, 2010.

\bibitem{JPS09}
G.~Jeronimo, D.~Perrucci, and J.~Sabia.
\newblock On sign conditions over real multivariate polynomials.
\newblock {\em Discrete Comput. Geom.}, 44(1):195--222, 2010.

\bibitem{GPT2013}
G.~Jeronimo, D.~Perrucci, and E.~Tsigaridas.
\newblock On the minimum of a polynomial function on a basic closed
  semialgebraic set and applications.
\newblock {\em SIAM J. Optim.}, 23(1):241--255, 2013.

\bibitem{KMSS}
H.~Kaplan, J.~Matou\v{s}ek, M.~Sharir, and S.~Safernov\'a.
\newblock Unit distances in three dimensions.
\newblock {\em Combinat. Probab. Comput.}, 21:597--610, 2012.

\bibitem{Lasserre01}
J.~B. Lasserre.
\newblock Global optimization with polynomials and the problem of moments.
\newblock {\em SIAM J. Optim.}, 11(3):796--817 (electronic), 2000/01.

\bibitem{Lasserre06-2}
J.~B. Lasserre.
\newblock Convergent {SDP}-relaxations in polynomial optimization with
  sparsity.
\newblock {\em SIAM J. Optim.}, 17(3):822--843 (electronic), 2006.

\bibitem{Lasserre06-1}
J.~B. Lasserre.
\newblock A sum of squares approximation of nonnegative polynomials.
\newblock {\em SIAM J. Optim.}, 16(3):751--765 (electronic), 2006.

\bibitem{Laurent09}
M.~Laurent.
\newblock Sums of squares, moment matrices and optimization over polynomials.
\newblock In {\em Emerging applications of algebraic geometry}, volume 149 of
  {\em IMA Vol. Math. Appl.}, pages 157--270. Springer, New York, 2009.

\bibitem{Loj2}
S.~Lojasiewicz.
\newblock Triangulation of semi-analytic sets.
\newblock {\em Ann. Scuola Norm. Sup. Pisa, Sci. Fis. Mat.}, 18(3):449--474,
  1964.

\bibitem{Milnor2}
J.~Milnor.
\newblock On the {B}etti numbers of real varieties.
\newblock {\em Proc. Amer. Math. Soc.}, 15:275--280, 1964.

\bibitem{monk}
L.~G. Monk.
\newblock {\em Elementary-recursive decision procedures}.
\newblock PhD thesis, UC Berkeley, 1975.

\bibitem{NN93}
Y.~Nesterov and A.~Nemirovskii.
\newblock {\em Interior-point polynomial algorithms in convex programming},
  volume~13 of {\em SIAM Studies in Applied Mathematics}.
\newblock Society for Industrial and Applied Mathematics (SIAM), Philadelphia,
  PA, 1994.

\bibitem{Pap}
C.~Papadimitriou.
\newblock {\em Computational Complexity}.
\newblock Addison-Wesley, 1994.

\bibitem{Parrilo03}
P.~A. Parrilo.
\newblock Semidefinite programming relaxations for semialgebraic problems.
\newblock {\em Math. Program.}, 96(2, Ser. B):293--320, 2003.
\newblock Algebraic and geometric methods in discrete optimization.

\bibitem{Perrucci09}
D.~Perrucci.
\newblock Linear solving for sign determination.
\newblock {\em Theoret. Comput. Sci.}, 412(35):4715--4720, 2011.

\bibitem{OP}
I.~G. Petrovski{\u\i} and O.~A. Ole{\u\i}nik.
\newblock On the topology of real algebraic surfaces.
\newblock {\em Izvestiya Akad. Nauk SSSR. Ser. Mat.}, 13:389--402, 1949.

\bibitem{KP97}
L.~Porkolab and L.~Khachiyan.
\newblock On the complexity of semidefinite programs.
\newblock {\em J. Global Optim.}, 10(4):351--365, 1997.

\bibitem{Powers-Reznick-2001}
V.~Powers and B.~Reznick.
\newblock A new bound for {P}\'olya's theorem with applications to polynomials
  positive on polyhedra.
\newblock {\em J. Pure Appl. Algebra}, 164(1-2):221--229, 2001.
\newblock Effective methods in algebraic geometry (Bath, 2000).

\bibitem{Powers-Wormann98}
V.~Powers and T.~W{\"o}rmann.
\newblock An algorithm for sums of squares of real polynomials.
\newblock {\em J. Pure Appl. Algebra}, 127(1):99--104, 1998.

\bibitem{Ramana97}
M.~V. Ramana.
\newblock An exact duality theory for semidefinite programming and its
  complexity implications.
\newblock {\em Math. Programming}, 77(2, Ser. B):129--162, 1997.
\newblock Semidefinite programming.

\bibitem{Rannou98}
E.~Rannou.
\newblock The complexity of stratification computation.
\newblock {\em Discrete Comput. Geom.}, 19(1):47--78, 1998.

\bibitem{Reif79}
J.~Reif.
\newblock Complexity of the mover's problem and generalizations.
\newblock In {\em IEEE Transactions on Robotics and Automation}, pages
  421--427, 1979.

\bibitem{R92}
J.~Renegar.
\newblock On the computational complexity and geometry of the first-order
  theory of the reals. {I}-{I}{I}{I}.
\newblock {\em J. Symbolic Comput.}, 13(3):255--352, 1992.

\bibitem{Rotman}
J.~J. Rotman.
\newblock {\em An Introduction to Algebraic Topology}.
\newblock Springer-Verlag, 1988.

\bibitem{RS}
M.-F. Roy and A.~Szpirglas.
\newblock Complexity of the computations with real algebraic numbers.
\newblock {\em Journal of Symbolic computation}, 10:39--51, 1990.

\bibitem{MS03}
M.~Safey El~Din and {\'E}.~Schost.
\newblock Polar varieties and computation of one point in each connected
  component of a smooth algebraic set.
\newblock In {\em Proceedings of the 2003 {I}nternational {S}ymposium on
  {S}ymbolic and {A}lgebraic {C}omputation}, pages 224--231 (electronic), New
  York, 2003. ACM.

\bibitem{MS04}
M.~Safey El~Din and {\'E}.~Schost.
\newblock Properness defects of projections and computation of at least one
  point in each connected component of a real algebraic set.
\newblock {\em Discrete Comput. Geom.}, 32(3):417--430, 2004.

\bibitem{Mohab-Schost2010}
M.~Safey El~Din and {\'E}.~Schost.
\newblock A baby steps/giant steps probabilistic algorithm for computing
  roadmaps in smooth bounded real hypersurface.
\newblock {\em Discrete Comput. Geom.}, 45(1):181--220, 2010.

\bibitem{Mohab-Schost2014}
M.~Safey El~Din and {\'E}.~Schost.
\newblock A nearly optimal algorithm for deciding connectivity queries in
  smooth and bounded real algebraic sets.
\newblock {\em CoRR}, abs/1307.7836, 2013.

\bibitem{SS}
J.~Schwartz and M.~Sharir.
\newblock On the piano movers' problem ii. general techniques for computing
  topological properties of real algebraic manifolds.
\newblock {\em Adv. Appl. Math.}, 4:298--351, 1983.

\bibitem{Solymosi-Tao}
J.~Solymosi and T.~Tao.
\newblock An incidence theorem in higher dimensions.
\newblock {\em arXiv:1103.2926v2 [math.CO]}.

\bibitem{Tarasov-Vyalyi08}
S.~P. Tarasov and M.~N. Vyalyi.
\newblock Semidefinite programming and arithmetic circuit evaluation.
\newblock {\em Discrete Appl. Math.}, 156(11):2070--2078, 2008.

\bibitem{Tarski51}
A.~Tarski.
\newblock {\em A decision method for elementary algebra and geometry}.
\newblock University of California Press, Berkeley and Los Angeles, Calif.,
  1951.
\newblock 2nd ed.

\bibitem{T}
R.~Thom.
\newblock Sur l'homologie des vari\'et\'es alg\'ebriques r\'eelles.
\newblock In {\em Differential and Combinatorial Topology (A Symposium in Honor
  of Marston Morse)}, pages 255--265. Princeton Univ. Press, Princeton, N.J.,
  1965.

\bibitem{Toda}
S.~Toda.
\newblock P{P} is as hard as the polynomial-time hierarchy.
\newblock {\em SIAM J. Comput.}, 20(5):865--877, 1991.

\bibitem{Wuthrich76}
H.~R. W\"uthrich.
\newblock {Ein Entschiedungsverfahren f\"ur die Theorie der
  reell-abgeschlossenen K\"orper}.
\newblock {\em LNCS}, 43:138--162, 1976.

\bibitem{Zahl}
J.~Zahl.
\newblock An improved bound on the number of point-surface incidences in three
  dimensions.
\newblock {\em Contrib. Discrete Math.}, 8:100--121, 2013.

\end{thebibliography}
\end{document}